
\magnification=\magstep1
\vsize=22truecm
\input amstex
\centerline{\bf{BALL CHARACTERIZATIONS IN PLANES}}
\centerline{\bf{AND SPACES OF CONSTANT CURVATURE, I}}
\vskip1cm
\centerline{J. JER\'ONIMO-CASTRO*, E. MAKAI, JR.**}
\vskip.3cm
\centerline{This pdf-file is not identical with the printed
paper.}
\vskip.3cm

*Research (partially) supported by CONACYT, SNI 38848
\newline
**Research (partially) supported by Hungarian National Foundation for 
Scientific Research, grant nos. T046846, T043520, K68398, K81146, Research
supported by ERC Advanced Grant
``GeoScape'', No. 882971.
\vskip1cm
ABSTRACT.
High proved the following theorem. 
If the intersections of any two congruent
copies of a plane convex body are centrally symmetric, then
this body is a
circle. 
In our paper we extend the theorem of High to the sphere and
the hyperbolic plane, and partly to 
spaces of constant curvature. We also investigate the dual
question about the
convex hull of the unions, rather than the intersections.

Let us have in $S^2$, ${\bold{R}}^2$ 
or $H^2$
a pair of convex bodies (for $S^2$ different from $S^2$),
such that 
the intersections of any congruent copies of them are
centrally symmetric. Then our bodies are congruent circles.
If the intersections of any congruent copies of them are
axially symmetric, then our bodies are (incongruent) circles.

Let us have in $S^2$, ${\bold{R}}^2$  
or $H^2$ proper
closed convex subsets $K,L$ with interior points,
such that 
the numbers of the connected
components of the boundaries of $K$ and $L$ are finite.
If the intersections of any congruent copies of $K$ and $L$ are centrally 
symmetric, then $K$ and $L$
are congruent circles, or, for ${\bold{R}}^2$, 
parallel
strips.
We exactly describe all 
pairs of such subsets $K,L$, whose any
congruent copies have an intersection with axial symmetry
(for $S^2$, ${\bold{R}}^2$ 
and $H^2$
there are 1, 5 and 9 cases, resp.).

Let us have in $S^d$, ${\bold{R}}^d$ 
or $H^d$ proper
closed convex $C^2_+$ subsets $K,L$ with interior points,
such that all sufficiently small intersections of their
congruent copies are symmetric w.r.t.\ a particular
hyperplane. Then the boundary components of both $K$ and $L$
are congruent, and each of them is
a sphere, a parasphere or a hypersphere.

Let us have a pair of
convex bodies in $S^d$, ${\bold{R}}^d$ 
or $H^d$,
which have at any boundary points supporting
spheres (for $S^d$ of
radius less than $\pi /2$). If
the convex hull of the union of any congruent copies of
these bodies is
centrally symmetric, then our bodies are congruent balls
(for $S^d$
of radius less than $\pi /2$). An analogous statement holds
for symmetry w.r.t.\ a particular hyperplane.
For $d=2$, suppose the
existence of the above supporting circles
(for $S^2$ of radius less than $\pi /2$),
and, for $S^2$, smoothness of $K$ and $L$.
If we suppose axial symmetry of all
the above convex hulls, then our bodies are (incongruent)
circles
(for $S^2$ of radii less than $\pi /2$). (This abstract
concerns both parts I and II of this paper.)

\vskip.1cm

$^*$ 
Facultad de Ingenier\'\i a, Universidad Aut\'onoma de Quer\'etaro, Centro
Uni\-ver\-si\-ta\-rio, 
\newline Cerro de las Campanas s/n C.P. 76010,
Santiago de Quer\'etaro, Qro. M\'exico, ME\-XI\-CO
\newline
{\rm{ORCID ID: https://orcid.org/0000-0002-6601-0004}}
\newline
$^{**}$ Alfr\'ed R\'enyi Institute of Mathematics,
Hungarian Research Network (HUN-REN),
\newline
H-1364 Budapest, Pf. 127, HUNGARY
\newline
{\rm{http://www.renyi.hu/\~{}makai}}
\newline
{\rm{ORCID ID: https://orcid.org/0000-0002-1423-8613}}

\vskip.1cm

{\it{E-mail address:}}
$^*$ jeronimo\@cimat.mx, jesusjero\@hotmail.com
\newline
$^{**}$ makai.endre\@renyi.hu

\vskip.5cm

1. INTRODUCTION

\vskip.1cm

We write
$S^d$, ${\bold{R}}^d$, 
$H^d$, with $d \ge 2$,
for the $d$-dimensional spherical,
Euclidean and hyperbolic spaces, resp.
{\it{Convexity of a set $K \subset H^d$}} is defined
as for $K \subset {\bold{R}}^d$. 
{\it{Convexity of $K \subset S^d$, with}} 
int\,$K \ne \emptyset $, is meant as follows:
for
any two non-antipodal points of $K$ the shorter great circle
arc 
connecting them belongs to
$K$. Then for $\pm x \in K$, $y \in $ int\,$K$ and
$y \ne \pm x$,
the shorter arcs ${\widehat{(\pm x)y}}$
belong to $K$, hence some half large circle connects $\pm x$
in $K$.
By a {\it{convex body}} 
in $S^d$, ${\bold{R}}^d$, 
$H^d$ 
we mean a compact convex set, with nonempty
interior. 
{\it{In}} $S^d$, when saying {\it{ball, or sphere}}, 
we always mean one {\it{of radius at most}} $\pi /2$ (thus a
ball is convex).
A proper closed convex subset
of $S^d$, ${\bold{R}}^d$  
or $H^d$, with nonempty interior,
is {\it{strictly
convex}}, if its boundary does not contain a non-trivial
segment. A {\it{convex surface}} is the boundary of a proper
closed convex subset 
of $S^d$, ${\bold{R}}^d$  
or $H^d$ with nonempty interior. For
$d = 2$ a convex surface will be called a {\it{convex curve}}.


R. High proved the following theorem.


\vskip.1cm

{\bf{Theorem.}} ([7])
{\it{Let $K \subset {\bold{R}}^2$  
be a convex body. Then
the following statements are equivalent:
{\rm{(1)}} All intersections $(\varphi K) \cap (\psi K)$, 
having interior points,
where $\varphi, \psi :
{\bold{R}}^2  
\to {\bold{R}}^2$  
are congruences, are centrally
symmetric.
{\rm{(2)}} $K$ is a circle.}} {\bf{QED}} 

\vskip.1cm

It seems, that his proof gives the analogous statement, when $\varphi , \psi $
are only allowed to be orientation preserving congruences.


\vskip.1cm

{\bf{Problem 1.}} 
Describe the pairs of closed convex sets with interior points, in $S^d$,
${\bold{R}}^d$,  
and $H^d$, whose any congruent copies have a centrally symmetric
intersection, provided this intersection has interior points.
Evidently, two 
congruent balls (for $S^d$ of radii at most $\pi /2$), or two parallel slabs
in ${\bold{R}}^d$, 
have a centrally symmetric intersection, provided it has a
nonempty interior.

\vskip.1cm


It was proved in [8],
Theorem 2, that in the $C^2$
case for $S^d$, and in the
$C^2_+$ case for ${\bold{R}}^d$ 
and $H^d$,
the
only possibility is two
congruent balls (for $S^d$ of radii at most $\pi /2$).
(For the case of ${\bold{R}}^d$, 
the solution cf.\ in the
end of this introduction.)


The authors are indebted to L. Montejano (Mexico City) and G. Weiss (Dresden)
for having turned their interest to characterizations of pairs of convex
bodies with all translated/congruent copies having a centrally or axially
symmetric intersection or convex hull of the union, resp., 
or with other symmetry properties, e.g., having some affine symmetry.


The aim of our paper is to give partial answers to this
problem. To exclude trivialities, we always
suppose that {\it{our sets are different from the whole plane, or space}},
and also we investigate only such cases, when the 
{\it{intersection has interior points}}.
We prove the analogue of the theorem of High for $S^2$ and $H^2$.
Namely, 
we characterize the pairs of proper
closed convex subsets with interior points, in $S^2$,
${\bold{R}}^2$ 
and 
$H^2$, having
centrally symmetric intersections of all congruent
copies, provided these intersections have nonempty
interiors.
However, for $H^2$ we have to suppose that if the connected  
components of the boundaries of
both subsets are straight lines, then there are altogether
finitely many of them.
Also we investigate a variant of this question, for $S^2$,
${\bold{R}}^2$ 
and $H^2$,
when we prescribe not central but 
axial symmetry of all intersections, having nonempty
interiors. Nevertheless, for $H^2$
we have to suppose that
for both subsets the number of connected components
of the boundary is finite.
We exactly describe all pairs of proper closed convex
subsets with
interior points, with the above property: for $S^2$,
${\bold{R}}^2$ 
and $H^2$
there are $1$, $5$ and $9$ cases, resp.

Suppose that in $S^2$, ${\bold{R}}^2$ 
and
$H^2$, all small 
intersections of congruent copies of two closed convex
proper subsets with
interior points, having a nonempty interior, admit some
non-trivial congruence. Then all connected components of the boundaries of the
two sets are cycles or straight lines.

The dual question is the question of central symmetry of
the closed convex hull 
of any congruent copies of $K$ and $L$. Under the
hypotheses that both $K$ and $L$ have at any of their boundary
points supporting spheres (for $S^d$ of radii less than $ \pi /2 $), 
the only case is two congruent balls (for $S^d$ of radii less than $\pi /2$).
This slightly sharpens Theorem 4 of [8].
(Observe that for $S^d$, ${\bold{R}}^d$ 
and $H^d$ this hypothesis implies that
any existing sectional curvature of $K$ and $L$ is positive, positive, or
greater than $1$, in the three cases, resp.)
Restricting ourselves to the case $d = 2$, prescribe
axial symmetry of all these closed convex hulls. Suppose 
that $K$ and $L$ have at all boundary
points supporting circles (for $S^2$ of
radius less than $\pi /2$), and for $S^2$ both $K$ and $L$
are $C^1$. Then the only case is two (incongruent) circles
(for $S^2$ of radii smaller than $\pi /2$).

Surveys about characterizations of central symmetry, for convex
bodies in
${\bold{R}}^d$, 
cf.\ in [3],
\S 14, pp. 124-127, and, more recently, in [6],
\S 4.

In a paper under preparation [10]
we will give more detailed theorems about
${\bold{R}}^d$. 
We will describe the pairs of proper closed convex subsets with interior points, without regularity hypotheses,
whose any congruent copies have (1) a centrally symmetric intersection (provided
this intersection has interior points); (2)
a centrally symmetric closed convex hull of their union.
For $d \ge 2$ in case (1) we have congruent
balls, or (incongruent) parallel slabs. For 
$d \ge 2$ in case (2) we have either infinite cylinders over 
balls of dimensions $2 \le i,j \le d$ with equal radii
(e.g., congruent balls),
or one subset is an infinite cylinder with dimension of axis
$0 \le i \le d - 1$ and with base compact, and the other
subset is a slab (e.g., (incongruent) parallel slabs).
These results form additions to the results of the papers [14], [15].

This introduction concerns both parts I and II of this paper.
Part I contains 1. Introduction, 2. New results: Theorems
1--4, 3. Preliminaries, 4. Proofs of Theorems 1--4.
Part II (i.e., [9])
contains Abstract (about Part II), 5. New results:
Theorems 5--8, 6. Preliminaries (about Part II), 7. Proofs of
Theorems 5--8.

\vskip.5cm

2. NEW RESULTS, THEOREMS 1--4

\vskip.1cm

We mean by a {\it{non-trivial congruence}} a congruence
different from the identity. 
We write conv\,$(\cdot )$, diam\,$(\cdot )$,
${\text{int}}\,( \cdot )$, ${\text{cl}}\,( \cdot )$,
${\text{bd}}\,( \cdot )$ and ${\text{perim}}\,( \cdot )$
for the convex hull, diameter,
interior, closure, boundary and perimeter
of a set.
A {\it{paracircle}} (also called a horocircle) is a
closed convex set in $H^2$, bounded by a paracycle.

As a general hypothesis in our theorems about intersections,
we have that 
$$
X {\text{ is }} S^d, \,\,{\bold{R}}^d 
{\text{ or }} H^d,
{\text{ with }} d \ge 2, 
{\text{ and }} K,L \subset 
X {\text{ are closed
convex proper subsets}}
$$
$$
{\text{with interior points. Moreover, }}
\varphi , \psi : X \to X, {\text{ sometimes with indices,}}
$$
$$
{\text{are orientation preserving congruences, with }}
{\text{int}}\,[ (\varphi K) \cap (\psi L) ]
\ne \emptyset .
\tag 1
$$
%

Sometimes we will say direct/indirect congruence for 
orientation preserving/re\-vers\-ing congruence.

The following Theorem
is the basis of our considerations for the case of
intersections, for $d=2$. 

\vskip.1cm


{\bf{Theorem 1.}}
{\it{Assume (1)
with $d = 2$.
Then we have $(1) \Longrightarrow (2)
\Longrightarrow (3)$, where
%

{\rm{(1)}}
There exists some $\varepsilon (K,L) > 0$, such that 
for each $\varphi , \psi $, for 
which
\,\,\,{\rm{diam}}\,
\newline
$[ (\varphi K) \cap (\psi L) ]
\le \varepsilon (K,L)$, 
we have that $(\varphi K) \cap (\psi L)$ is axially
symmetric.
%

{\rm{(2)}}
There exists some $\varepsilon (K,L) > 0$, such that 
for each $\varphi , \psi $, for 
which\,\,\,{\rm{diam}}\,
\newline
$[ (\varphi K) \cap (\psi L) ]
\le \varepsilon (K,L)$, 
we have that $(\varphi K) \cap (\psi L)$
admits some non-trivial congruence.
%

{\rm{(3)}}
Each connected component of the boundaries of both $K$ and
$L$ is a
cycle (for $X=S^2$ a circle of radius at most $\pi /2$),
or a straight line. If either for $K$, or for $L$, one
connected component is a circle, or paracycle, then this is
the unique component, and $K$, or $L$ is a circle (disk), or a
paracircle, resp.  
%

In particular, if the congruences in {\rm{(2)}} are central symmetries,
then in {\rm{(3)}} the 
connected components of the boundaries of both $K$ and $L$
are congruent.
For $X=S^2$ and $X = {\bold{R}}^2$ 
we have
$(1) \Longleftrightarrow (2) \Longleftrightarrow (3)$.
For $X=H^2$, if both for $K$ and $L$, the infimum of the positive
curvatures of its boundary components 
is positive, and at most one of its boundary components has
$0$ curvature, 
then $(1) \Longleftrightarrow (2) \Longleftrightarrow (3)$.

Let $X=H^2$. 
If for, e.g., $K$, 
the infimum of the positive curvatures
of its boundary components is $0$, or 
two of its boundary components have $0$ curvatures, then
$(3) \not\Longrightarrow (2)$.
Even,
supposing $(3)$ for $K$, we may prescribe in any way
the curvatures of the connected hypercycle or straight line
boundary
components of $K$ (with multiplicity), in case that
the infimum of the positive curvatures of the boundary
components of $K$
is $0$, or two boundary
components of $K$ have
$0$ curvatures. 
Then we can find a $K$ with these prescribed curvatures of the
connected hypercycle or straight line boundary
components of $K$ (with multiplicity), and
an $L$, such that for them $(3)$ holds, but $(2)$
does not hold.}} 


\vskip.1cm

As follows from Theorem 1,
the compact case is particularly simple. This of course
includes the case when $X = S^2$.

\vskip.1cm


{\bf{Theorem 2.}}
{\it{Assume {\rm{(1)}}
with $d = 2$.
Let both $K$ and $L$ be compact. Alternatively, as
a particular case of this, let $X = S^2$.

Then we have $(1) \Longleftrightarrow (2) \Longleftrightarrow
(3)$, where
%

{\rm{(1)}}
For each $\varphi , \psi $
we have that $(\varphi K) \cap (\psi L)$
is centrally symmetric.
%

{\rm{(2)}}
There exists some $\varepsilon (K,L) > 0$, such that for each 
$\varphi , \psi $, for which 
{\rm{diam}}\,
\newline
$[ (\varphi K) \cap (\psi L) ] \le
\varepsilon (K,L)$,
we have that
$(\varphi K) \cap (\psi L)$ is centrally symmetric.
%

{\rm{(3)}}
$K$ and $L$
are congruent circles, for $S^2$  
of radius at most $\pi /2$.
%
Also we have $(4) \Longleftrightarrow (5) \Longleftrightarrow
(6) \Longleftrightarrow (7) \Longleftrightarrow  (8)$, where
%

{\rm{(4)}}
For each $\varphi , \psi $
we have that $(\varphi K) \cap (\psi L)$
is axially symmetric.
%

{\rm{(5)}}
For each $\varphi , \psi $
we have that $(\varphi K) \cap (\psi L)$
admits some non-trivial congruence.
%

{\rm{(6)}}
There exists some $\varepsilon (K,L) > 0$, such that
for each $\varphi , \psi $, for which
\,\,\,{\rm{diam}}\,
\newline
$[ (\varphi K) \cap (\psi L) ]
\le \varepsilon (K,L)$, 
we have that $(\varphi K) \cap (\psi L)$
is axially symmetric.
%

{\rm{(7)}}
There exists some $\varepsilon (K,L) > 0$, such that
for each $\varphi , \psi $, for 
which\,\,\,{\rm{diam}}\,
\newline
$[ (\varphi K) \cap (\psi L) ]
\le \varepsilon (K,L)$, 
we have that $(\varphi K) \cap (\psi L)$
admits some non-trivial congruence.
%

{\rm{(8)}}
$K$ and $L$ are {\rm{(in general incongruent)}} circles, for $S^2$
of radii at most $\pi /2$.}}
%


\vskip.1cm

{\bf{Theorem 3.}}
{\it{Assume (1)
with $d = 2$ and let $X = {\bold{R}}^2$. 
Then we have $(1) \Longleftrightarrow (2)$, where

{\rm{(1)}}
For each $\varphi , \psi $
we have that $(\varphi K) \cap (\psi L)$
admits some non-trivial congruence.

{\rm{(2)}}
$K$ and $L$ are {\rm{(in general incongruent)}}
circles, or one of them is a circle and
the other one is a
parallel strip or a 
half-plane, or they are two parallel strips, or they are two
half-planes.

In particular, writing in {\rm{(1)}}
central
symmetries (rather than non-trivial congruences)
is equivalent to writing in {\rm{(2)}} either two
congruent circles or two parallel strips.
Similarly, writing in {\rm{(1)}} axial symmetries is
equivalent to adding to {\rm{(2)}}
that for the case of two parallel
strips, these strips are congruent.}}


\vskip.1cm

The following two theorems give two different characterizations
for $H^2$, 
under different additional hypotheses. Of these, Theorem 4
deals with central symmetry, and Theorem 5 (in Part II, i.e., [9])
deals with axial
symmetries, or non-trivial congruences.
Recall Theorem 1,
$(2) \Longrightarrow (3)$.


\vskip.1cm

{\bf{Theorem 4.}}
{\it{Assume (1)
with $d = 2$ and let $X=H^2$. 
If all connected components of the boundaries of both of $K$
and $L$ are
straight lines, let their total number be finite.  
Then we have $(1) \Longleftrightarrow (2)$, where 

{\rm{(1)}}
For each $\varphi , \psi $
we have that $(\varphi K) \cap (\psi L)$
is centrally symmetric.

{\rm{(2)}}
$K$ and $L$ are congruent circles.}}

\vskip.1cm


{\bf{Problem 2.}}
Is the finiteness hypothesis in Theorem 4 necessary?


\vskip.1cm

{\bf{Remark 1.}} 
As will be seen from the proof of Theorem 4, namely in the
proof of Lemma 4.1, rather than the finiteness hypothesis
in Theorem 4, we may suppose only the following. There holds
one of the following properties.
\newline
(1) 
One of $K$ and $L$, e.g.,
$K$ has a boundary component $K_1$, with the
following
property. Let us pass
on ${\text{bd}}\, K$, meant in $B^2$ containing the model
circle, in the positive sense.
Then there is a non-empty open arc $A_K$ of $S^1$, which
immediately follows $K_1$ on this boundary
(thus $K_1$ and $A_K$
have a common infinite point), and which
does not contain any infinite point of any
boundary component of $K$. Simultaneously, 
the other
set $L$ has a boundary component $L_1$, with the following
property. Let us pass
on ${\text{bd}}\, L$, meant in $B^2$ containing the model
circle,
in the negative sense.
Then there is a non-empty open arc $A_L$ of $S^1$, which immediately
follows $L_1$ on this boundary (thus $L_1$ and $A_L$
have a common infinite point), and which
does not contain any infinite point of any
boundary component of $L$.
\newline
(2) The same as (1), only with ``positive sense'' and 
``negative sense'' interchanged.
\newline
As will follow from Part II (i.e., [9]),
Theorem 5, namely from the proofs of
Lemmas 5.7--5.9, rather than the finiteness hypothesis in
Theorem 4, we may suppose the following.
\newline
(3) Either $K$, or $L$ has two boundary components with at
least one common infinite point. (Cf.\ Part II (i.e., [9]),
Remark 3.)


\vskip.1cm

In the proofs of our Theorems we will use some ideas of [7].


\vskip.5cm

3. PRELIMINARIES

\vskip.1cm

We write $\| \cdot \| $ for the norm of a vector, and $B^d$
for the closed unit ball, in ${\bold{R}}^d$. 
For $x,y$
in $S^d, {\bold{R}}^d$ 
or $H^d$, we write $d(x,y)$ for their
distance in the respective space, and
$[x,y]$ or $(x,y)$ for the closed or open
segment (shorter segment in $S^d$) with end-points $x,y$. (We
will not apply this last notation for antipodal points on
$S^d$.) For $x$ a point, and $A$ a subset of $S^d,
{\bold{R}}^d$ 
or $H^d$, we write dist\,$(x,A)$ for the distance
of $x$ to $A$.
The line $xy$ is the line spanned by $x,y$ (this
notation will not be applied for $x = y$, or for $S^d$ and
$x + y = 0$). Suppose for $d = 2$, that
$x_1,x_2$ on the boundary of a closed
convex set $K \subset X$ with interior points are
``close'' to each other. Then
we write
${\widehat{x_1x_2}}$ {\it{for the (shorter, or unique)
arc of ${\text{\rm{bd}}}\,K$, with these endpoints}}, which
set $K$ will be clear from the context. Also we will
specify, whether we mean a closed or an open arc.

We recall that if the boundary of
a closed convex set $K \subset {\bold{R}}^d$ 
with interior points
is differentiable, then it is $C^1$. By the collinear models
of $S^d$ and $H^d$, this takes over to $S^d$ and $H^d$. We
will say in this case that $K$ {\it{is}} $C^1$.

We write $X$ 
for $S^d$, ${\bold{R}}^d$, 
$H^d$ for $d \ge 2$. (Except in
Theorems 6 and 7, we will be concerned with the case $d = 2$.)
For $x \in X$ 
and $r > 0$, we write $B(x,r) \subset X$ for the
{\it{closed ball in $X$,
of centre $x$ and radius $r$}} (for $X = S^d$ supposing $r \le
\pi /2$). For $Z \subset Y \subset X$
we write ${\text{relint}}_Y Z$ for the relative interior of
$Z$ w.r.t.\ $Y$.

For hyperbolic plane geometry we refer to [2], [4], [11], [12].
For geometry
of hyperbolic
space we refer to [1] 
[5] 
For elementary differential geometry we
refer to [16].

We shortly recall some of the concepts to be used later.

Two straight lines in $H^2$ can have either at least one common
finite point, then they are {\it{intersecting}}; or at least
one 
common infinite point, then they are {\it{parallel}} (in
these cases coincidence is allowed); or none
of
these cases occurs, then they are {\it{ultraparallel}}.

In
$S^2$, ${\bold{R}}^2$ 
and
$H^2$ there are the following (maximal, connected, twice 
differentiable) curves of constant curvature (in $S^2$ meaning geodesic
curvature). In $S^2$ these are 
the circles, of radii $r \in (0, \pi /2]$, with (geodesic)
curvature $\cot r
\in [0, \infty )$.
In ${\bold{R}}^2$, 
these are 
circles of radii $r \in (0, \infty )$, with curvature $1/r \in
(0, \infty )$; and straight lines, with curvature $0$.
In $H^2$, these are circles of radii $r \in (0, \infty )$,
with curvature coth$\,r \in (1, \infty )$; paracycles, with curvature $1$; and
hypercycles, i.e., distance lines, with distance $l \in (0,
\infty )$ from their base lines
(i.e., the straight lines
that connect their points at infinity), with curvature tanh$\,l \in (0,1)$;
and straight lines, with curvature $0$. 
Paracycles are also called horocycles, and their unique
infinite points are also called their {\it{centres}}.
Either in $S^2$, ${\bold{R}}^2$, 
or in $H^2$,
each sort of the above curves have
different curvatures, and for one sort, with different $r$ or $l$, they also
have different curvatures.
The common name of these curves 
is, except for straight lines in ${\bold{R}}^2$ 
and $H^2$, 
{\it{cycles}}. In $S^2$ also a great
circle is called a {\it{cycle}}, but when speaking about straight lines, for
$S^2$ this will mean great circles. 
An elementary method for the calculation
of these curvatures in $H^2$ cf.\ in [17].

Sometimes we will include straight lines among the
hypercycles (with $l= 0$).
Then the base line of a straight line is meant to be itself.
In this case cycles occur as orthogonal trajectories of all
straight lines incident to a finite point $p$ (all circles with
centre $p$); of all
straight lines incident to an infinite point $q$ (all
paracycles
with infinite point $q$); of all straight lines orthogonal to
a straight line $l$ (all hypercycles with base line $l$).

The space $H^d$ has two usual models, in
${\text{int}}\, B^d$,
namely the
collinear (Beltrami-Caley-Klein)
model (cf.\ [18]),
and the conformal
(Beltrami-Poincar\'e) model (cf.\ [20]).
(Sometimes we will consider the closed unit
ball, when we have to consider the infinite points as well.
Then we say closure of the model in ${\bold{R}}^d.)$ 
In analogy, we will speak about collinear and conformal
models of $S^d$ in ${\bold{R}}^d$. 
By this we mean the ones obtained by central
projection (from the centre), 
or by stereographic projection (from the north pole), to the
tangent hyperplane of $S^d$,
at the south pole, in ${\bold{R}}^{d+1}$. 
These exist of
course only on the open southern half-sphere, or on $S^d$ minus the north pole,
resp. Their images are ${\bold{R}}^d$. 
We  call the
{\it{centre of the model}} the south pole of $S^d$.
The collinear and conformal models of ${\bold{R}}^d$ 
are meant as
itself, with {\it{centre}} the origin.

Sometimes we will consider the (collinear or conformal) model
circle of $H^2$ as the unit circle of the complex plane
${\bold{C}}$. 
Thus
we will speak about its points $1$, $i$, etc.

A {\it{paraball}} (also called a horoball)
is a closed convex set in $H^d$,
bounded by a parasphere.

The congruences of $S^2$, ${\bold{R}}^2$ 
and $H^2$ can be given as follows
(cf.\ [1],
p. 70, and Theorems 4.1 and 4.2, where $H^2$ denotes the Poincar\'e
upper halfplane model, [1],
p.\ 43).
The orientation preserving ones are rotations in $S^2$, rotations and
translations in ${\bold{R}}^2$, 
and rotations, ``rotations about an infinite
point'', and translations along a straight line (preserving this line) in 
$H^2$. The
orientation reversing ones are glide reflections in each of 
$S^2$, ${\bold{R}}^2$ 
and $H^2$. For $H^2$ each congruence can be
uniquely extended by continuity to the closure of the
(collinear or conformal) model circle in
${\bold{R}}^2$, 
to a homeomorphism of this closure. 

For a topological space $Y$ we say that some property of a
point
$y \in Y$ holds {\it{generically}}, if it holds outside a
nowhere dense closed subset.

\vskip.5cm

4. PROOFS OF THEOREMS 1--4

\vskip.1cm

In the proofs of our theorems by the 
{\it{boundary components of a set}} we will
mean the connected components of the boundary of that set.

The following (2), (3), (4)
and
(5)
will serve to exclude in $H^2$ in many cases
congruences of $(\varphi K) \cap
(\psi L)$, which are non-trivial translations, or 
glide reflections which are not reflections, or  
non-trivial rotations about an infinite point, or rotations
about a finite point which are not central symmetries.

If a closed convex proper subset $M \subset  H^2$
with nonempty interior admits a glide reflection
which is not a reflection, as a congruence to itself,
then it admits the square of this glide reflection as well,
which is a non-trivial translation. (Therefore we will not
need to exclude glide reflections which are
not reflections, but exclusion of non-trivial translations
will suffice.)
If $M$ admits a non-trivial translation,
then it contains the closed convex hull of
the orbit of some point, w.r.t.\ the subgroup generated by
this translation.
Thus, $M$ contains a straight line. Consequently
$$
{\text{Let a closed convex proper subset $M \subset H^2$
with
nonempty interior admit a non-trivial}}
$$
$$
{\text{translation, or a glide reflection which is not a
reflection. Then $M$ has two}}
$$
$$
{\text{different infinite points (in particular, it
is not compact).}}
\tag 2 
$$
In an analogous way we have
$$
{\text{Let a closed convex proper subset $M \subset H^2$
with
nonempty interior admit a non-trivial}} 
$$
$$
{\text{rotation about an infinite point. Then $M$ contains a
paracircle with this}}
$$
$$
{\text{infinite point (centre; in particular, $M$
is not compact).}}
\tag 3 
$$

Now let $x_1x_2x_3x_4$ be a strictly convex quadrangle in 
$H^2$. Suppose that it admits a congruence $\chi $ which is a
{\it{combinatorial central
symmetry}} of our quadrangle
(i.e., $\chi x_1 = x_3$, $\chi x_2 = x_4$, 
$\chi x_3 = x_1$ and $\chi x_4 = x_2$). Then $\chi $ is
orientation preserving.
Further, $\{ x_1,x_3 \} $ is invariant under $\chi $, hence
also $[x_1,x_3]$ is invariant, and also its midpoint $o$ is
invariant. Since $\chi $ is orientation preserving, it is a
rotation about $o$, and from above the angle of rotation is
$ \pi $.
Hence $\chi $ is a central symmetry
of our quadrangle,
w.r.t.\ a unique centre $o$.
In fact, if there were two such
centres, then the
composition of the two central symmetries
would be a non-trivial translation admitted by our quadrangle,
contradicting
(2).
Rephrasing this,
$$
{\text{Let $x_1x_2x_3x_4$ be a strictly convex quadrangle in 
$H^2$. Suppose that}}
$$
$$
{\text{it admits a congruence $\chi $ which is combinatorial
central symmetry}}
$$
$$
{\text{(cf.\ above). Then $\chi $ is a central symmetry
w.r.t.\ a unique centre $o$.}}
\tag 4 
$$

Further, let a closed convex proper subset $M \subset H^2$
with
nonempty interior 
have a connected boundary component $M'$ of its
boundary, with the following property. It is
invariant either under a
non-trivial translation, or a glide reflection which is not a
reflection, or a non-trivial rotation about an infinite point.
Suppose that $M'$ has a nonsmooth point. Then it has
infinitely many
nonsmooth points. Rephrasing this,
$$
{\text{Let a closed convex set $M \subset H^2$
with
nonempty interior have a boundary}}
$$
$$
{\text{component $M'$, which has a positive finite number of
nonsmooth points.}}
$$
$$
{\text{Then a non-trivial congruence admitted by $M'$ is
either a rotation about}}
$$
$$
{\text{a finite point (e.g., a central symmetry), or is a
reflection w.r.t.\ a straight}}
$$
$$
{\text{line (according to as the congruence is orientation
preserving, or orientation}}
$$
$$
{\text{reversing).}} 
\tag 5 
$$

Let $X = S^2$, ${\bold{R}}^2$ 
or $H^2$.
Let $K \subset X$ be a closed convex set with interior
points. Let
$x \in {\text{bd}}\,K$, and let $x'$
be a point of ${\text{bd}}\,K$ very close to $x$, that
follows $x$ on bd$\,K$
in the positive sense.
We will often consider the {\it{shorter
counterclockwise arc}}
${\widehat{xx'}}$ 
{\it{of}} ${\text{bd}}\,K$. This makes sense if $K$ is
compact, i.e., ${\text{bd}}\,K$ is
homeomorphic to $S^1$. If
the connected  component of the boundary of $K$, containing
$x$, 
is homeomorphic to ${\bold{R}}$, 
and thus connects two
possibly coincident
infinite points, then there is
just one such arc. Later, 
when writing {\it{shorter arc}}, we mean the shorter
one in the first case (the other arc has a length almost
the perimeter of $K$) and the unique one in the second
case.

The {\it{distortion of the arc element in $X = S^d$,
${\bold{R}}^d$ 
or $H^d$, in the
collinear, or conformal model, resp.}}, is the quotient of
the corresponding
arc element in the collinear, or conformal model, resp.,
as a subset of ${\bold{R}}^d$, 
and of the arc element in
$X$. For $x,y \in X$ we write $x',y'$ for
their images in the collinear or conformal
model, and $d'(x',y')$ for
the distance of $x',y'$ in the collinear, or conformal model,
resp., as a subset of ${\bold{R}}^d$. 
(We will always tell,
which model do we mean.)

We recall ([18])
and ([20])
that in the collinear, or the conformal model of $H^d$, resp.,
the arc element at $x' \in {\text{int}}\, B^d$ is given by
$$
ds^2 = \| dx' \| ^2 / (1 - \| x' \| ^2 ) +
(\langle x', dx' \rangle )^2/
(1 - \| x' \| ^2 )^2 \in
\left[ \| dx' \| ^2 , \| dx' \| ^2 \times \right.
$$
$$
\left. [ 1/(1 - \| x' \| ^2) + \| x' \| ^2/
(1 - \| x' \| ^2)^2 ] \right] , {\text{ or }} ds^2 = 4 
\| dx' \| ^2 / (1 - \| x' \| ^2 ) ^2, {\text{ resp.}}
$$
$$
{\text{So, in compact sets }} C \subset H^d,
{\text{ the distortion of the arc element in }} H^d,
{\text{ in the}}
$$
$$
{\text{collinear, or conformal model, resp.,
is bounded below and above. Hence,}}
$$
$$
{\text{for distinct }}
x,y \in C, {\text{ we have that }} d'(x',y')/d(x,y) {\text{ is
bounded below and}}
$$
$$
{\text{above.}}
\tag 6 
$$

The first statement of the following
(7)
is elementary.
$$
{\text{In compact sets }} C {\text{ of the open southern
hemisphere of }} S^d, {\text{ the}}
$$
$$
{\text{distortion of the arc element in }} S^d, {\text{ in the
collinear, or conformal}}
$$
$$
{\text{model, resp., is bounded below and above. Hence, for
distinct }} x,y
$$
$$
\in C, {\text{ we have that }} d'(x',y')/d(x,y) {\text{ is
bounded below and above.}}
\tag 7 
$$

Now we prove the last statements in
(6)
and in
(7).

Observe that in the collinear model geodesic segments
(shorter geodesic segments in the open southern hemisphere
of $S^d$) are
preserved, hence the distances, in $X$, or in ${\bold{R}}^d$, 
can be obtained by integrating
the respective arc-elements along them. This proves the
last statements for the collinear model.

For the conformal model we may suppose that $C \subset X$ is
a ball, whose centre is mapped to $0$ in the model, and
for the case of $S^d$ that
its radius is smaller than $\pi / 2$.
Thus $C \subset X$ is convex, and also its image $C'$ in the
conformal model is convex (it is a ball of centre $0$).
Now let us consider for $x,y \in C$
the (shorter) geodesic
segment $[x,y]$ in $X$, which lies in $C$. Its image in the
model is a curve
joining $x'$ and $y'$, having a length at most
${\text{const}}_C \cdot d(x,y)$, hence $d'(x',y') \le
{\text{const}}_C \cdot d(x,y)$.
Changing the role of $C$ and $C'$,
in the same way we gain $d(x,y) \le
{\text{const}}_C \cdot d'(x',y')$.


\vskip.1cm

{\it{Proof of Theorem {\rm{1}}.}} 
{\bf{1.}}
The implication $(1) \Longrightarrow (2)$ in Theorem 1 is
evident.

We turn to the proof of the implication $(2)
\Longrightarrow (3)$ in Theorem 1. This will be finished by
Lemma 1.7.

\vskip.1cm


The following lemma is surely known, but we could not
locate a proof for it. Therefore we give its simple proof.
It is some analogue of
(6)
and
(7).


\vskip.1cm

{\bf{Lemma 1.1.}}
{\it{Let $X = S^2$ or $X = H^2$.
Let $o \in X$ be fixed, and let its image in the collinear
model be the centre of the model. Let $p \in X$, let $d(o,p)
\le r$,  and let us
consider an
angle in $X$ with apex $p$. Then the quotient of the
measure of the image of our angle
in the collinear model, as a subset of ${\bold{R}}^2$, 
and of the
measure of this angle in $X$
lies in $[\cos r, 1/\cos r]$ for $X = S^2$ and
$r < \pi /2$, and in
$[1/{\text{\rm{cosh}}}\,r, {\text{\rm{cosh}}}\,r]$ for
$X = H^2$. In both cases the lower
and upper bounds are sharp.}}


\vskip.1cm

{\it{Proof.}}
{\bf{1.}}
We begin with the case of $S^2$. Its collinear model is
obtained by central projection in ${\bold{R}}^3$  
of the open
southern hemisphere, from the origin, to the tangent
plane of $S^2$ at the south pole. We may assume  $d(o,p)
= r$, and $p = (\sin r, 0, - \cos r)$. In the tangent plane
of $S^2$ at $p$ we consider an orthonormal coordinate system
with basis vectors $e_1' := (\cos r, 0, \sin r )$ and $e_2' :=
(0, 1, 0)$. We consider a rotating unit vector in this
coordinate system, given by $f(\Phi ) :=
(e_1' \cos \Phi, e_2' \sin \Phi )$, for
$\Phi \in [0, 2 \pi ]$. Let $\varepsilon \in (0, \pi /2 -
r)$.
Then consider $p$ and
$$
p + \tan \varepsilon \cdot f(\Phi ) = (\sin r +
\tan \varepsilon \cdot \cos r \cdot \cos \Phi ,
\tan \varepsilon \cdot \sin  \Phi ,
- \cos r + \tan \varepsilon \cdot
\sin r \cdot \cos \Phi )
\tag 8 
$$
(which lies in the open half-space $z < 0$ by
$\varepsilon < \pi /2 - r$)
and their images $p'$ and $[ p + \tan \varepsilon \cdot
f(\Phi )] '$
by the map $(x, y, z) \mapsto (- x/z, -y/z, -1)$
from the open lower half-space to the plane $z = -1$.
(This map, when
restricted to the open southern hemisphere of $S^2$,
gives the map to the collinear model. The image of
$p + \tan \varepsilon \cdot f(\Phi )$ by this map
runs over the image of the small
circular line
on $S^2$, of centre $p$, and radius $\varepsilon $.)
Then 
$p' = (\tan r, 0, - 1)$, and
$$
[ p + \tan \varepsilon \cdot f(\Phi ) ] ' =
\left( \frac{\sin r + \tan
\varepsilon \cdot \cos r \cdot \cos \Phi }
{\cos r - \tan \varepsilon \cdot \sin r \cdot \cos \Phi },
\frac{\tan \varepsilon \cdot \sin \Phi }
{\cos r - \tan
\varepsilon \cdot \sin r \cdot \cos \Phi }, -1 \right) .
\tag 9 
$$
Substracting from the last expression $p'$, after some
simplification we get
$$
\left(
\frac{\tan \varepsilon \cdot \cos \Phi }{(\cos r - \tan
\varepsilon \cdot
\sin r \cdot \cos \Phi \cdot \cos r) \cos r},
\frac{\tan \varepsilon \cdot \sin \Phi }{\cos r - \tan
\varepsilon \cdot
\sin r \cdot \cos \Phi \cdot \cos r},
0
\right) .
\tag 10 
$$
Then the slope of the line with this direction vector in the
plane is
$\tan \Phi \cdot \cos r$. The angle of this line with the basic
vector $e_1'$ is $\arctan (\tan \Phi \cdot \cos r)$.

Obviously it is sufficient to prove the statement of the
lemma
for ``infinitesimally small'' angles. Then integration will
prove the same inequality for ``finite'' angles.
Therefore, differentiating $\arctan (\tan \Phi
\cdot \cos r)$
w.r.t.\ $\Phi $, after some simplifications this derivative
becomes
$$
(\cos r) /(1 - \sin ^2 r \cdot \sin ^2 \Phi ) \in 
[\cos r, 1/ \cos r] .
\tag 11 
$$
Here both endpoints of the interval
are attained, for $\Phi = 0, \pi $,
and for $\Phi = \pi /2, 3 \pi /2$, resp. 
(The $\tan $ function has a singularity at 
$\Phi = \pi /2, 3 \pi /2$. However, then we can exchange
the
basic vectors $e_1'$ and $e_2'$, and the slope changes to
its reciprocal. Then we can use the
$\cot $ function of the angles in $X$ and in the model,
which is analytic there. Thus if we choose the branches of the
$\arctan $ function so that $\arctan (\tan t \cdot \cos r)$
remains continuous at $(k + 1/2) \pi $ for $k$ any integer,
then $\arctan (\tan t \cdot \cos r)$
becomes a strictly increasing
analytic function on ${\bold{R}}$.) 

\vskip.1cm

{\bf{2.}}
We turn to the case of $H^2$.
We may suppose $d(o,p) = r$ and
$p = ({\text{tanh}}\, r, 0)$ in polar coordinates,
by the definition
of the collinear model of $H^2$. Consider a circle of centre
$p$ and of radius $\varepsilon > 0$.
Again we take a rotating segment $[p,q]$, with $d(p,q) =
\varepsilon $,
from the centre $p$ of this circle to its boundary point $q$,
with $\angle opq = \pi - \Phi $, where $\Phi
\in [0, 2 \pi ]$.
We denote by $'$ the
points, quantities in the collinear model corresponding to
points, quantities in $H^2$.

We consider the positively oriented triangle $opq$.
We denote $\Psi := \angle poq \in (0, \pi )$
and $s := d(o,q)$. Consider
$r,s, \Psi $ as given. Then we have by the cosine law
${\text{cosh}}\, \varepsilon = {\text{cosh}}\, r \cdot
{\text{cosh}}\, s - {\text{sinh}}\, r \cdot {\text{sinh}}\, s
\cdot \cos \Psi $, and by the sine law $\sin \Phi =
{\text{sinh}}\, s \cdot 
(\sin \Psi )/ {\text{sinh}}\, \varepsilon
= {\text{sinh}}\, s \cdot (\sin \Psi )
/ ({\text{cosh}}^2\, \varepsilon - 1)^{1/2}$,
where we substitute the value of ${\text{cosh}}\, \varepsilon
$ from above. Last we calculate $|\tan \Phi | = \sin
\Phi /
(1 - \sin ^2 \Phi)^{1/2}$.

In the collinear model we have the image $0p'q'$ of the
triangle $opq$. Here $\Psi ' = \Psi $ and $r' =
{\text{tanh}}\, r$ and $s' = {\text{tanh}}\, s$. Like above,
we determine first the side $\varepsilon '$ of our triangle,
then $\sin \Phi ' $, and last
$|\tan \Phi '|$. We claim
$$
\tan \Phi ' = \tan \Phi \cdot {\text{cosh}}\, r.
\tag 12 
$$

Since the signs of $\tan \Phi $ and $ \tan \Phi '$ are
the same, it suffices to show the squared equality. To show
this, we perform the calculations indicated above, expressing
everything with the variables $r,s,\Psi $. We 
cancel
with $\sin \Psi = \sin \Psi '$, and clear all the
denominators. Thus we obtain two equal quantities, quadratic
polynomials of $\cos \Psi $, with coefficients depending on
$r$ and $s$, namely
$$
{\text{sinh}}^2\, r \cdot {\text{cosh}}^2\, s - 2 
{\text{cosh}}\, r \cdot {\text{sinh}}\, r \cdot
{\text{cosh}}\, s \cdot {\text{sinh}}\, s \cdot \cos \Psi +
{\text{cosh}}^2\, r \cdot {\text{sinh}}^2\, s \cdot 
\cos ^2 \Psi .
\tag 13 
$$
Thus
(12)
is proved.

Again, we need to calculate 
the derivative $d\Psi ' /d\Psi = (d/d\Psi ) \arctan (\tan
\Psi \cdot {\text{cosh}}\, r)$, and to determine
its minimum and maximum. This derivative is
$$
({\text{cosh}}\, r)/(1 + {\text{sinh}}^2\, r \cdot \sin ^2
\Psi ) \in [1/{\text{cosh}}\, r, {\text{cosh}}\, r] .
\tag 14 
$$
Here the left endpoint of the interval is attained for
$\Psi = \pi /2$, and the right endpoint is asymptotically
attained for $\Psi \to 0$ and for $\Psi \to \pi $.
{\bf{QED}}

\vskip.1cm


{\bf{Lemma 1.2.}}
{\it{Assume \thetag{1} 
with $d = 2$.
Then {\rm{(2)}} of Theorem {\rm{1}} implies that
both $K$ and $L$ are $C^1$.}}


\vskip.1cm

{\it{Proof.}} 
{\bf{1.}}
A closed convex set $K \subset
S^2, {\bold{R}}^2, 
H^2$ 
with interior points is not necessarily
differentiable. However, at each boundary point $x$
it has two 
half-tangents, that is, the limit of the line $xx'$, when
$x' \to x$,  so that $x'$ follows $x$ on ${\text{bd}}\,K$
in the positive, or in the negative sense. 

Suppose,
e.g., that $K$ is not $C^1$. Let $x \in {\text{bd}}\,K$ be
a point of
non-smoothness. Let $\alpha \in (0, \pi )$ 
denote the angle of the positively oriented half-tangents
of $K$ at $x$. We will call this angle the {\it{outer angle
of $K$ at $x$}}. (This is $\pi $ minus the inner angle
of $K$ at $x$.)

Let $x', x'' \in {\text{bd}}\,K$
be points very close to $x$, such that the shorter, say,
counterclockwise
open arc $\widehat{x'x''}$ contains $x$. 
Furthermore,
we choose the points $x',x''$ so that, additionally, 
we have $d(x,x') = b \varepsilon $ and 
$d(x,x'') = c \varepsilon $. Here $\varepsilon > 0$ is
small, and
$b,c \in (0, 1]$ fixed satisfy, that 
a Euclidean triangle $T$ 
with one angle $\pi - \alpha $ and adjacent sides $b,c$ is
not isosceles.
Observe that, by continuity, all sufficiently small
distances occur as $d(x,x')$ and $d(x,x'')$. 
Similarly, let $y,y' \in {\text{bd}}\,L$ be such that $y'$
is very close to $y$, and $y'$ follows $y$ on ${\text{bd}}\,L$
in the positive
sense. Like above, all sufficiently small
distances occur as $d(y,y')$.
Let $d(x',x'') = d(y,y')$.

Then there exist orientation preserving congruences $\varphi $
and $\psi $, with the following properties. We have
$\varphi (x') = \psi (y')$, and
$\varphi (x'') = \psi (y)$, and 
$(\varphi K) \cap (\psi L)$ is bounded by the shorter
arcs ${\widehat{\varphi (x') \psi (x'')}}$ 
of ${\text{bd}}\,(\varphi K)$ and ${\widehat{\psi (y) \psi
(y')}}$ of ${\text{bd}}\,(\psi L)$. 
Thus this intersection is an arc-triangle $A$,
with ``vertices'' $\varphi x , \varphi (x') =
\psi (y')$ and $\varphi (x'') = \psi (y)$.
Let $T$ be the triangle
with the same vertices.

We are going to prove that
$$
{\text{any congruence admitted by }} A {\text{ preserves the
set of its three ``vertices''.}} 
\tag 15 
$$

\vskip.1cm

{\bf{2.}}
First we deal with the case $X = {\bold{R}}^2$.  
By one-sided
differentiablity of ${\text{bd}}\,K$ at $x$, and of 
${\text{bd}}\,L$ at $y$, the total angular rotations (i.e.,
curvature measures) of all
the three open arc-sides of $A$ are $o(1)$, for $\varepsilon
\to 0$.
In particular, the half-tangents of $A$ at its ``vertices''
enclose with the respective
half-tangents of $T$ at the same vertices,
i.e., with the
respective side lines of $T$, an angle $o(1)$. (I.e., the
half-tangents enclose small angles with the secant lines,
whose limits are the half-tangents.)
Therefore the inner
angles of $A$ differ from the respective inner angles of $T$
by $o(1)$. Thus the inner angles of $T$ are $\pi - \alpha
+ o(1)$, and two other angles, whose sum is $\alpha
+ o(1)$. Therefore both of these last mentioned
inner angles of $T$ are at most $\alpha
+ o(1)$, and the respective outer angles are therefore at
least $\pi - \alpha + o(1)$. Therefore each outer angle of $A$
at its ``vertices'' is at least $\min \{ \pi - \alpha + o(1),
\alpha + o(1) \} $.

Now let us consider one of the open arc-sides of $A$.
Consider all the points $p$ in this open arc-side of
$A$, such that the outer angle at $p$ is positive. The sum
of all the outer angles at such points $p$ is at most the
total angular rotation of the considered open arc-side,
which is $o(1)$ for $\varepsilon \to 0$. Hence for each open
arc-side of $A$, all the points $p$ in this open arc-side of
$A$, for which the outer angle is positive, 
satisfy that this outer angle is $o(1)$, for
$\varepsilon \to 0$.

Therefore, for $\varepsilon > 0$ sufficiently small, the
outer angles of $A$ are either large, namely at least 
$\min \{ \pi - \alpha + o(1), \alpha + o(1)\} $, or small,
namely $o(1)$. Therefore any congruence admitted by
$A$ preserves the
three largest outer angles of $A$, which occur at the
``vertices'' of $A$. This proves
(15)
for $X =
{\bold{R}}^2$. 

\vskip.1cm

{\bf{3.}}
Second we deal with the case of $X = S^2, H^2$. Here total
angular rotation makes no sense, therefore we have to go to
${\bold{R}}^2$ 
via the collinear model. Denote the images of
$A$ and $T$ in the collinear model by $A'$ and $T'$, resp.
{\it{Let $\varphi x$ be
mapped to the centre $0$ of the collinear
model.}} Then, for $\varepsilon
> 0$ sufficiently small, the total angular rotations in
${\bold{R}}^2$ 
of the open arc-sides of $A'$ are arbitrarily
small. Like in {\bf{2}}, this implies that in
${\bold{R}}^2$, 
the half-tangents of $A'$ at its ``vertices''
enclose with the respective
half-tangents of $T'$ at the same vertices,
i.e., with the
respective side lines of $T'$, an angle $o(1)$. 
Therefore, like in {\bf{2}}, in ${\bold{R}}^2$, 
the inner
angles of $A'$ differ from the respective inner angles of $T'$
by $o(1)$.
This implies, like in {\bf{2}}, that in ${\bold{R}}^2$, 
each outer angle of $A'$
at its ``vertices'' is at least $\min \{ \pi - \alpha + o(1),
\alpha + o(1) \} $.

However, the collinear model distorts the angles.
By Lemma 1.1,
we have lower and upper bounds of the form
$1 + o(1)$ for the ratio of the angles in 
$X$ and their images in ${\bold{R}}^2$, 
for $\varepsilon \to 0$.
Therefore the outer angles of $A$ at its ``vertices'' in $X$
are at least $(\min \{ \pi - \alpha + o(1),
\alpha + o(1) \} ) \cdot [ 1 + o(1) ] =
\min \{ \pi - \alpha + o(1), \alpha + o(1) \} $.

Now let us consider one of the open arc-sides of $A$.
Consider all the points $p$ in this open arc-side of
$A$, such that the outer angle at $p$ is positive (this means
the same for $X$ and for the images in ${\bold{R}}^2$). 
The sum
of all the outer angles, at the images in ${\bold{R}}^2$ 
of such points $p$, 
is at most the
total angular rotation of the considered open arc-side 
in ${\bold{R}}^2$, 
which is $o(1)$ for $\varepsilon \to 0$. Hence for each open
arc-side of $A$, all the points $p$ in this open arc-side of
$A$, for which the outer angle is positive, 
satisfy that the image of
this outer angle in ${\bold{R}}^2$ 
is $o(1)$, for
$\varepsilon \to 0$. Then the same angle in $X$ is at most
$o(1) \cdot [ 1 + o(1) ] = o(1)$. Hence, like in
{\bf{2}}, 
any congruence admitted by $A$ preserves the
three largest outer angles of $A$, which occur at the
``vertices'' of $A$. This proves
(15)
for
$X = S^2, H^2$.

\vskip.1cm

{\bf{4.}}
We turn once again to $X = {\bold{R}}^2$ 
The inner
angle of $T$ at $\varphi x $ is $\pi - \alpha + o(1)$.
If the side of $T$ opposite to this angle is
$a \varepsilon $, then
$
a^2 = b^2 + c^2 - 2bc
\cos [ \pi - \alpha + o(1) ] = b^2 + c^2 - 2bc
\cos ( \pi - \alpha ) + o(1)$.

Observe that ${\text{int}}\, A \ne \emptyset $, since it has
an inner angle in $(0, \pi )$. 
We claim that also
diam\,$A$ is small. Clearly diam\,$A$ is
attained for a pair of points on the arc-sides of $A$.
For $T \,\,(\subset A)$ we have that
diam\,$T$ is at most $\max \{ a,b,c \} \cdot \varepsilon $,
which is small. Now it suffices to observe that an arc-side
of $A$ has a distance $o(\varepsilon )$
from the respective side of $T$. This
follows from the fact that the angles of the arc-sides of $A$
and the respective sides of $T$ at both of their endpoints
are $o(1)$, hence the
arc-sides are contained in isosceles triangles with base
the respective side of $T$, and height $o(1)$. Therefore
diam\,$A \le [ \max \{ a,b,c \} + o(1) ]
\varepsilon $, thus is arbitrarily small. 
Since also ${\text{int}}\, A \ne \emptyset $, therefore $A$
admits a non-trivial congruence.

If for some sequence of
$\varepsilon $'s, tending to $0$, the arc-triangle $A$ admitted
a
non-trivial congruence, then it would preserve its
``vertices''.
Hence it
would be a non-trivial congruence admitted by $T$ as well. That is,
$T$ would be an isosceles triangle.
Then the limit triangle, satisfying $a^2 = 
b^2 + c^2 - 2bc \cos ( \pi - \alpha )$, would be isosceles
too, contradicting the choice of $b$ and $c$.

Hence for all
sufficiently small $\varepsilon $ the triangle $T$ is not
isosceles, and thus $A$ admits no non-trivial congruence,
contradicting the hypothesis of Lemma 1.2 (i.e.,
(2) of Theorem 1).

\vskip.1cm

{\bf{5.}}
We turn once again to
$X = S^2, H^2$. Analogously as in {\bf{4}}, now we
have to write the spherical and hyperbolic
cosine laws for the side of length 
$a \varepsilon $ of the triangle $T$ in $X$.
Like usual, we subtract
$1$ from both sides of this equation, and then divide
both sides by
$\varepsilon ^2$. Thus we obtain 
$a^2 = b^2 + c^2 - 2bc
\cos [ \pi - \alpha + o(1) ] + O( \varepsilon ^2)
= b^2 + c^2 - 2bc \cos ( \pi - \alpha ) + o(1)$.
This is an analogous formula as in the beginning of {\bf{4}}.

As in {\bf{4}}, ${\text{int}}\,A \ne \emptyset $. Now we have
diam\,$T \le \max \{ b, c \} \cdot 2 \varepsilon $,
since $T$ can be included in a circle of centre $\varphi x$
and radius $\max \{ b, c \} \cdot \varepsilon $. Also
now the
arc-sides of $A$ have a distance at most $o(\varepsilon )$
from the respective sides of $T$. This follows by including
the arc-side to an isosceles triangle like in {\bf{4}},  
and using spherical and hyperbolic trigonometric formulas.
Then the diameter of $A$ can occur between two points of $T$,
or one point of $T$ and one point in the ``half-lens like
domains'' between
the sides of $T$ and the respective arc-sides of $A$, or
between two points in such ``half-lens like
domains''. In each case we have
diam\,$A \le O(\varepsilon ) + o(\varepsilon ) =
o(1)$.

The limit argument is the same as for
${\bold{R}}^2$, 
hence
we obtain a contradiction once more.

{\bf{6.}}
The conclusions of {\bf{4}} and {\bf{5}} contradict our
indirect hypothesis about non-smooth\-ness of $K$. This proves
smoothness of $K$, and of $L$.
{\bf{QED}}


\vskip.1cm

{\bf{Lemma 1.3.}}
{\it{Assume
{\rm{(1)}}
with $d = 2$. Let $K$ be $C^1$. For $K$
compact let $K'' := K' := {\text{\rm{bd}}}\, K$ {\rm{(this is
homeomorphic to $S^1$)}}. For $K$ non-compact let $K'$ be a
connected component of ${\text{\rm{bd}}}\,K$ {\rm{(this is
homeomorphic to ${\bold{R}}$, 
and tends to infinity at both
of its ``ends'', for ${\bold{R}}^2$ 
and $H^2$)}},
moreover, let $K''$ be a
compact subarc of $K'$. Then there exists an $\varepsilon
(K,K'') > 0$ such that for each $x \in
K''$ and each $\varepsilon \in (0,
\varepsilon (K,K'')]$ there hold the following statements.
\newline
{\rm{(1)}}
The intersection $({\text{\rm{bd}}}\,K) \cap
B(x, \varepsilon )$ is a closed subarc of
$K'$.
\newline
{\rm{(2)}}
The intersection $({\text{\rm{bd}}}\,K)
\cap {\text{\rm{bd}}}\, [
B(x, \varepsilon ) ]$ consists of exactly two points
$x_{\varepsilon }^+$ and $x_{\varepsilon }^-$, the
endpoints of the subarc in {\rm{(1)}}. The
directed arcs of
${\text{\rm{bd}}}\,K$ from $x$ to $x_{\varepsilon }^+$ and $x_{\varepsilon }^-$ are
positively and negatively directed, resp.
\newline
{\rm{(3)}}
The angles of ${\text{\rm{bd}}}\,K$ and ${\text{\rm{bd}}}\, [
B(x, \varepsilon ) ] $,
at the two points of intersection from {\rm{(2)}}, are 
$\pi /2 + o(1)$, for $\varepsilon \to 0$, uniformly for each
$x \in K''$.
\newline
{\rm{(4)}}
The functions $x \mapsto x_{\varepsilon }^+$ and 
$x \mapsto x_{\varepsilon }^-$ are continuous, and
are strictly monotonous in the following sense. A small
motion of $x$ along $K''$ strictly in positive or
negative sense implies
small motions of $x_{\varepsilon }^+$ and
$x_{\varepsilon }^-$ along $K''$ strictly in positive or
negative sense, resp.}}


\vskip.1cm

{\it{Proof.}}
{\bf{1.}}
We consider $K$ and $K''$ as fixed. Hence dependence of
quantities on $K$ and $K''$ will
not be explicitly written in the notations.

We begin with the case when ${\text{bd}}\,K$
{\it{is connected}}
(i.e., ${\text{bd}}\,K = K'$).

{\it{Let $\delta > 0$ be sufficiently small.}} We define
$K''(\delta ) \subset K'$ as follows. {\it{For $K$
compact}}, we let $K''(\delta ) = {\text{bd}}\, K$.
{\it{For $K$ noncompact}} first we extend $K''$ to the closed 
$2 \delta $-neighbourhood $K''_{2\delta }$ of $K''$ in $K'$,
and then to the closed
$4 \delta $-neighbourhood $K''_{4\delta }$ of $K''$ in $K'$,
in the arc length metric of $K'$.
Then we extend $K''_{4\delta }$
further as follows.
Observe that now $K'$ is homeomorphic to ${\bold{R}}$, 
and 
tends to infinity at both of its ``ends''. 
We take such a long compact subarc $K''(\delta )$
of $K'$, containing $K''_{4\delta }$,
that
we have
${\text{dist}}\,[ K''_{4\delta }, K' \setminus
K''(\delta ) ] \ge 1 $, with distance
meant in $X$.
This will ensure for $x \in K''$ and
$\delta < 1$ that
$$
K' \cap B(x, \delta ) =
\left[ [K''(\delta ) \cup [ K' \setminus K''(\delta ) ] \right]
\cap B(x, \delta ) =
K''(\delta ) \cap B(x, \delta ).
\tag 16 
$$
Therefore it will suffice to show (1), (2) and (3) of the
lemma, with $K'$ replaced by $K''(\delta )$.

\vskip.1cm

{\bf{2.}}
In this proof
we will use the collinear model, with coordinates
$\xi ,\eta $.
This clearly works for ${\bold{R}}^2$ 
and $H^2$, however for $S^2$ this exists
only on the open southern hemisphere. Let $X = S^2$.

First suppose diam\,$K = \pi $. Then either $K$ is a
halfsphere, or a digon with angle in $(0, \pi )$. In the first of these cases 
the statement of the lemma is evident, for any $\varepsilon \in (0, \pi )$.
In the second of these cases $K$ is not $C^1$.

Second suppose diam\,$K < \pi $. Then by the first four
sentences of the second paragraph of the proof of Lemma 1.4
in [8],
applied for $d = 2$, we get that
$K$ lies in an open
hemisphere. Then we may suppose that this is the open southern
hemisphere, hence the collinear model exists on some
neighbourhood of $K$.

\vskip.1cm

{\bf{3.}}
By
(6),
(7)
and Lemma 1.1, for $S^d$
and $H^d$, we have the following.
Both the arc elements and the
angles have quotients bounded below and above on
$K''( \delta )$ (with
vertex of the angle on $K''(\delta )$),
when considered in $X$, and in
the collinear model, as a subspace of
${\bold{R}}^2$. 
For
${\bold{R}}^2$ 
these are obvious.

We claim that
$$
{\text{for a sufficiently short subarc of }}
K''(\delta ), {\text{ the chord length and the arc length}}
$$
$$
{\text{in }} X
{\text{ have a quotient }}
1 + o(1),
{\text{ uniformly, if the arc length tends to }} 0.
\tag 17 
$$
For sufficiently short
arcs of $K''(\delta )$ the tangent lines in the collinear
model, as a subspace of${\bold{R}}^2$,  
change very little, uniformly. We will identify
${\text{bd}}\,K$ and $K''(\delta )$ by their images
in the collinear model, resp.
We cover $K''(\delta )$ by
four subsets, according to as the
the tangent direction
in the positive sense belongs to the open angular intervals
$(- \pi /3, \pi /3)$,
or $(\pi /6, 5\pi /6)$, or $(2\pi /3, 4 \pi /3)$, or
$(7 \pi /6, 11\pi /6)$. Thus we obtain four
open subsets $I_1, \dots , I_4$
of $K''(\delta )$, covering $K''(\delta )$.
On $I_1$ or $I_3$ we have that ${\text{bd}}\,K$ can be
given by an
equation $\eta = f(\xi )$, with $f$ convex, or concave,
resp. On $I_2$ or $I_4$ we have that
${\text{bd}}\,K$
can be given by an equation $\xi = g(\eta )$, with $g$
concave, or convex, resp. We have $|f'(\xi)|,
|g'(\eta )| < {\sqrt{3}}$. Moreover,
the domains of definition of
these functions in all four cases are the respective
projections of the open subsets $I_i$ of $K''(\delta )$
to the $\xi $- or $\eta $-axis, resp.

Let us consider a sufficiently short arc of
$K''(\delta )$.
Its first endpoint, in the positive sense, can have a
positively directed tangent
direction lying in $[- \pi /4, \pi /4]$, $[\pi /4,$
\newline
$3 \pi /4]$, $[3 \pi /4, 5 \pi /4]$ or 
$[5 \pi /4, 7 \pi /4]$. Suppose the first case.
(The other three
cases are settled analogously. If we can find an $\varepsilon
(K,K'')$ for the first case, then analogously we can find
$\varepsilon (K,K'')$ for the other three cases as well. 
Then the minimum of these four values will satisfy the
statement of the lemma.)
Then our
sufficiently short arc of $K''(\delta )$ lies in $I_1$.
Moreover, on it we have $\eta = f(\xi )$, with $f$
convex, and $|f'(\xi)| < {\sqrt{3}}$. 
Let the endpoints
of our sufficiently
short arc (and chord) be $(\xi _1, \eta _1)$ and
$(\xi _2, \eta _2)$.
Then our arc is given as $\{ \left( \xi , f(\xi ) \right) \mid
\xi \in [\xi _1, \xi _2 ] \} $.

The length of this arc is
$$
\int _{\xi_1} ^{\xi _2 } \left[ g_{11} \left( \xi , f(\xi )
\right) + 2 g_{12} \left( \xi , f(\xi ) \right) f'(\xi ) +
g_{22} \left( \xi , f(\xi ) \right]
[ f'(\xi ) ] ^2 \right] ^{1/2} d \xi .
\tag 18 
$$
Here $g_{ij}$ is the metric tensor, in the $(\xi , \eta )$
coordinate system.
Moreover,
the chord length is the same expression, with $f(\xi )$
replaced by ${\overline f}(\xi )
:= \eta _1 + (\eta _2 - \eta _1) \cdot (\xi - \xi _1)/(\xi _2 -
\xi _1)$, and hence with $ (\overline f)'(\xi ) = (\eta _2 - \eta _1)/
(\xi _2 - \xi _1)$.
(Since chords in $X$ and in the collinear model, as a subset
of ${\bold{R}}^2$, 
coincide.)
Here, by the mean value theorem, also using
that $f$ is $C^1$, on the interval $[\xi _1, \xi _2]$ we have
$f'(\xi ) - (\overline f)'(\xi ) = o(1)$, uniformly,
if the arc length tends to $0$ (hence also $\xi _2 - \xi _1$
tends to $0$). Hence on
this interval also $f(\xi ) - {\overline f}(\xi ) =
o(\xi _2 - \xi _1) = o(1)$,
uniformly,
if the arc length tends to $0$.

This implies that the difference of
(18),
and the analogous
expression, which is obtained from
(18)
by replacing $f(\xi )$
in it by ${\overline f}(\xi )$,
is $o(\xi _2 - \xi _1)$. E.g., we
show $g_{22}\left( \xi , f(\xi ) \right) \times $
\newline
$[ f'(\xi ) ] ^2 -  
g_{22}\left( \xi , {\overline f}(\xi ) \right) [
({\overline f})'(\xi ) ] ^2 = o(1)$, uniformly for $\xi
\in [\xi _1, \xi _2]$. (The analogous estimates for the first
and second summands of the square of the
integrand in
(18)
are obtained analogously, but even
simpler.) We have
$$
g_{22}\left( \xi , f(\xi ) \right) [ f'(\xi )
] ^2 -  g_{22}\left( \xi , {\overline f}(\xi ) \right)
[ ({\overline f})'(\xi ) ] ^2 = [g_{22}\left( \xi ,
f(\xi ) \right) 
$$
$$
 - g_{22}\left( \xi , {\overline f}(\xi ) \right) ] \cdot
[ f'(\xi ) ] ^2 + g_{22}\left( \xi , {\overline f}
(\xi ) \right) \cdot [\left( f'(\xi ) + 
({\overline f})'(\xi ) \right) ] \times
$$
$$
[\left( f'(\xi ) 
- ({\overline f})'(\xi ) \right) ] = o(1) \cdot O(1) +
O(1) \cdot O(1) \cdot o(1) = o(1).
\tag 19 
$$
It remains to observe that the function $t \mapsto t^{1/2}$
is Lipschitz on any closed subinterval of $(0, \infty )$.
This we apply to
the square of the integrand in
(18),
which has a
positive lower (and upper) bound.
In fact,
the metric tensor can be estimated from below by some
const$_{K''} \cdot (d\xi ^2 + d\eta ^2) \ge {\text{const}}_
{K''} \cdot d\xi ^2$. Hence the square of the integrand
in
(18)
is at least const$_{K''}$.

The last lower estimate 
yields a lower bound for
(18),
namely
const$_{K''}^{1/2} \cdot (\xi _2 - \xi _1)$.
From above the difference of
(18)
and of the analogous
expression, with $f(\xi )$ replaced by ${\overline f}(\xi )$,
is $o(\xi _2 - \xi _1)$.
{\it{These imply our claim}}
(17).

\vskip.1cm

{\bf{4.}}
The $C^1$ curve $K''(\delta )$ has an inherited
Riemannian submanifold (with boundary)
structure, with the intrinsic metric. 
This is isometric either to some circular line in
${\bold{R}}^2$, 
with the metric the length of the not longer
connecting arc on the circular line,
or to some closed interval of ${\bold{R}}$,  
with the metric
inherited from ${\bold{R}}$ 
(according to whether $K$ is compact, or not). We denote
this metric by $d_{\text{arc}}$.
On the other hand, we have the metric
$d_X$ in $X$, which for $K''(\delta )$ is the
chord length in $X$.
Then $d_{\text{arc}} \ge d_X$,
hence the identical map of the compact metric space
$(K''(\delta ),d_{\text{arc}})$ to the compact metric space
$(K''(\delta ),d_X)$ is continuous. Hence 
this
bijective map is a homeomorphism. Hence its inverse map
$(K''(\delta ),d_X) \to (K''(\delta ),d_{\text{arc}})$
also is continuous, with compact metric domain, hence it is
even uniformly continuous. Let $U_{\text{arc}}(\delta ) :=
\{ (k_1,k_2) \in K''(\delta ) \times K''(\delta )
\mid d_{\text{arc}}
(k_1,k_2) \le \delta \} $ and $U_X(\delta ) :=
\{ (k_1,k_2) \in K''(\delta ) \times K''(\delta )
\mid d_X(k_1,k_2) \le \delta \} $.
Therefore,
for each
$\delta > 0$ there exists a $\gamma (\delta ) > 0$,
such that
$U_X [ \gamma (\delta ) ]
\subset U_{\text{arc}}(\delta )$.
Therefore, 
$$
{\text{for }} x \in K'' {\text{ and }} x' \in
K''(\delta ) {\text{ and}} 
$$
$$
d_X(x,x') \le \gamma (\delta ) {\text{ we have }}
d_{\text{arc}}(x,x') \le \delta .
\tag 20 
$$

$$
{\text{If }} K {\text{ is compact, we still assume }}
\delta \le ({\text{perim}}\,K)/8.
\tag 21 
$$
In this case,
for $d_{\text{arc}}(x,x') \le 4 \delta \le
({\text{perim}}\, K)/2$ we have that
$d_{\text{arc}}(x,x')$ equals the integral of
$ds$ on the positively or negatively
oriented arc $\widehat{xx'}$.

We will see (in
(26)
that 
$\varepsilon (K,K'')$ can be chosen as
$\min \{ \delta , \gamma (\delta ) \} $,
for some suitably small $\delta > 0$.

\vskip.1cm

{\bf{5.}}
Now consider a sufficiently short arc $\widehat{xx'}$
of $K''(\delta )$, with $x \in K''$, and with
$x'$ following $x$ on
$K''(\delta )$ in the positive sense,
with arc length
$\delta $, say. (The case when $x'$ follows $x$ on
$K''(\delta )$ in the negative sense, can be settled
analogously.)
We investigate
the direction, in the collinear model, 
as a subset of${\bold{R}}^2$,  
of the positively
oriented segment (chord)
$S'$ with endpoints $x,x'$, with orientation
inherited from
the positively oriented $K''(\delta )$.
Further, we investigate
the direction of the positively oriented
tangent line $( l(x')
) '$ of $K''(\delta )$ at $x'$,
in the collinear model, 
as a subset of ${\bold{R}}^2$. 
By the mean value theorem, and the
$C^1$ property of $K''(\delta )$,
we have the following.
The angle, in the
collinear model, as a subset of ${\bold{R}}^2$, 
of the direction of the
oriented segment $S'$, and the direction of the
oriented tangent line $( l(x') ) '$
is $o(1)$, uniformly, if the arc length $\delta $ tends
to
$0$. Observe that our oriented segment and oriented tangent
line in the model, as a subspace of ${\bold{R}}^2$, 
have  $x'$ as a common point.

Therefore {\it{the direction of the
corresponding oriented segment $S := [x,x']$
and oriented tangent
line $l(x')$ of $K''(\delta )$ at $x'$,
both taken in $X$, have
in $X$ an angle, at $x'$,
which is also $o(1)$, uniformly}}, by Lemma 1.1.
Now observe that $S$ is a radius of the circle
$B \left( x, d_X( x', x ) \right)$ in $X$,
hence is perpendicular to
${\text{bd}} \,[ B(x, d_X( x', x ) ] $
in $X$. Therefore
$$
l(x') {\text{ encloses in }} X {\text{ an angle }} \pi /2 +
o(1) {\text{ with }} {\text{bd}} \,[ B(x, \delta ) ]
$$
$$
{\text{at }} x', {\text{ with }} o(1) {\text{ uniform, for }}
x \in K'' {\text{ and }} \delta \to 0 .
\tag 22 
$$
Also taking into consideration
(16),
{\it{this would prove}} (3) {\it{of the lemma,
provided we already knew}} (2) {\it{of the lemma}}.

Let $x \in K''$.
We write $r_x$ for the distance from $x$
in $X$,
and $s_x$ for the
arc length distance of a point $x' \in K''$
from $x$ (i.e., the length of the not
longer connecting arc ${\widehat{xx'}} \subset  K''$).
(We have
$s_x(x) = r_x(x) = 0$.) For $K$ compact
we have $s_x \le ({\text{perim}}\,K)/2$;
therefore in this case we will consider only subarcs of 
${\text{bd}}\,K$ not longer than $({\text{perim}}\,K)/2$.
This ensures that {\it{the integral of
$ds$ on such arcs equals $s_x$}} (cf.\ the sentence following
(21)).
The analogue of this italicized statement for
the noncompact case is obvious (there is only one such arc).
Then
by the trigonometry of Euclidean, spherical and
hyperbolic triangles ([19]),
we have the following.
Formula
(22)
implies
$$
{\text{uniformly for all }} r {\text{ in an interval of the
form }} [0,r_0], {\text{ where }} r_0 {\text{ is}}
$$
$$
{\text{small, that }}
1 \le (ds_x/dr_x)(r) = 1 + o(1), {\text{ with }} o(1)
{\text{ uniform. Hence}}
$$
$$
r_x \le s_x = r_x \cdot [ 1 + o(1) ] ,
{\text{ with }} o(1) {\text{ uniform,
for }} x \in K'' {\text{ and}}
$$
$$
\delta \to 0;
{\text{ hence }} s_x \le 2r_x, {\text{ for }} x \in K''
{\text{ and }} \delta  {\text{ sufficiently small.}}
\tag 23 
$$
Then
$$
1 \ge dr_x/ds_x = 1 + o(1) \ge 1/2, {\text{ with }} o(1) 
{\text{ uniform, for }} \delta 
{\text{ sufficiently small.}} 
$$
$$
{\text{Therefore passing with }} x'
{\text{ away from }} x \in K'', {\text{ along an arc of }}
K''(\delta ), 
$$
$$
{\text{in the positive sense, of arc length at most }} 4 \delta
{\text{ from }} x, {\text{ the distance}}
$$
$$
d_X(x,x' ) {\text{ strictly increases, for a {\it{uniform
sufficiently small }} }} \delta > 0 .
\tag 24 
$$
{\it{Now we fix this value of $\delta $.}}

Here we can attain $d_X$-distance any
$\varepsilon \in (0, \delta]$,
for some $x' \in K''_{2\delta } \subset
K''(\delta )$.
In fact, for $x \in K''$
and $x' \in K''_{2\delta } \subset K''(\delta )$ with
$d_{\text{arc}}(x,x') = 2 \delta $ -- such an $x'$ exists, by the
definition of $K''_{2\delta }$ --
we have by
(23)
$d_X(x,x') =
2 \delta [ 1 + o(1) ] \ge \delta \ge \varepsilon $.
By the strictly increasing property of $d_X(x,x')$,
on the arc of $K''_{2\delta }$, consisting of the $(x')$'s,
following $x$ on $K''_{2\delta }$ in the positive sense and
satisfying $d_{\text{arc}}(x,x') \le 2 \delta $, we have the following.
The distance
$d_X(x,x')$
can be equal to $\varepsilon $ only for
one point $x'$; which point $x'$ in
fact exists, as pointed out above. By
(16),
(2) of the
lemma would be proved, if we proved it with
$K''(\delta )$ rather than $K'$. This in turn
would be proved,
provided we knew already
(1) of the lemma, also with $K''(\delta )$ rather than
$K'$. As observed in
(22),
this would prove
also (3) of the lemma.

{\bf{6.}}
First we will investigate the case of non-compact $K$. Then
$K'$ tends to infinity at both of its ends (for
${\bold{R}}^2$ 
and $H^2$).

From
(23),
on a positively oriented arc of
$K''(\delta )$ with starting point $x \in K''$,
of arc length
$2\delta $, 
both of $r_x$ and $s_x$ are strictly monotonically increasing
$C^1$
functions of each other: $r_x = r_x^+(s_x)$, and
$s_x = s_x^+(r_x)$.
Therefore,
$$
{\text{for }} x \in K'' {\text{ and for }} \varepsilon
\in (0, \delta ], {\text{ we have }} A^+_{\varepsilon } :=
\{ x' \in K''(\delta ) \mid x' {\text{ follows }} x
$$
$$
{\text{on }} K''(\delta ){\text{ in the positive sense}} \}
\cap B(x,\varepsilon ) = \{ x' \in K''(\delta ) \mid x'
{\text{ follows }} x {\text{ on}}
$$
$$
K''(\delta ) {\text{ in
the positive sense, and }} d_{\text{arc}}(x,x') = s_x^+
\left( d_X(x,x') \right) \le s_x^+  (\varepsilon ) =
\varepsilon \times
$$
$$
[ 1 + o(1) ] \le
2 \varepsilon \le 2 \delta \} , {\text{ by
(23),
with }}
o(\cdot ) {\text{ uniform, for }} x \in K'' {\text{ and }}
\delta \to 0.
\tag 25 
$$
Then $A^+_{\varepsilon } \subset K''(\delta )$
is a positively oriented closed arc
of $K''(\delta )$, with starting point $x$ and of arc length
$\varepsilon \cdot [ 1 + o(1) ] $,
with $o(\cdot )$ uniform. If on $A^+_{\varepsilon }$ a variable
point $x''$ moves from $x$ in the positive orientation, then
$d_X(x,x'')$ is a strictly increasing continuous function of
the position of $x''$. Hence it assumes each value till its
maximum value $\varepsilon $ just once.
We let
$$
0 < \varepsilon \le \varepsilon (\delta )
:= \min \{ \delta, \gamma (\delta ) \} .
\tag 26 
$$
By
(16),
(1) of the
lemma would be proved, if we proved it with
$x' \in K''(\delta )$ rather than $x' \in K'$.

By
(20),
$$
{\text{for }} x \in K'' {\text{ and }} x' \in K''(\delta )
{\text{ and }} d_{\text{arc}}(x,x') > \delta 
$$
$$
{\text{we have }} x' \not\in B(x, \gamma (\delta ) ) \supset
B(x, \varepsilon (\delta ) ) \supset
B(x, \varepsilon ) .
\tag 27 
$$

There remains the case when
$x \in K''$ and $x' \in K''(\delta )$, and
$d_{\text{arc}}(x,x') \le \delta $. 
By
(25)
and
(26)
we have
$$
\{ x' \in
K''(\delta ) \mid x' {\text{ follows }} x {\text{ on }}
K''(\delta ) {\text{ in the positive sense,}}
$$
$$
{\text{and }}
d_{\text{arc}}(x,x') \le \delta \} \cap
B(x, \varepsilon ) =
\big\{ x' \in
K''(\delta ) \mid x' {\text{ follows }} x {\text{ on}}
$$
$$
K''(\delta )
{\text{ in the positive sense, and }}
d_{\text{arc}}(x,x') \le \min \{ \delta , s_x^+( 
\varepsilon ) \} \big\} .
\tag 28 
$$
Hence, by
(26),
(27)
and
(28)
we get
$$
\{ x' \in
K''(\delta ) \mid x' {\text{ follows }} x {\text{ on }}
K''(\delta ) {\text{ in the positive sense}} \}
\cap B(x, \varepsilon )
$$
$$
=
[\{ x' \in
K''(\delta ) \mid x' {\text{ follows }} x {\text{ on }}
K''(\delta ) {\text{ in the positive sense, and}}
$$
$$
d_{\text{arc}}(x,x') \le \delta \} \cap
B(x, \varepsilon ) ]
\cup
[\{ x' \in
K''(\delta ) \mid x' {\text{ follows }} x {\text{ on }}
K''(\delta ) {\text{ in}}
$$
$$
{\text{the positive sense, and }}
d_{\text{arc}}(x,x') > \delta \} \cap
B(x, \varepsilon ) ] 
=
\big\{ x' \in
K''(\delta ) \mid x'
$$
$$
{\text{follows }} x {\text{ on }}
K''(\delta ) {\text{ in the positive sense, and }}
d_{\text{arc}}(x,x') \le \min \{ \delta , s_x^+( 
\varepsilon ) \} \big\} .
\tag 29 
$$

Let us turn to points $x'$ following $x$ on
$K''(\delta )$
in the negative sense. Then analogously to
(25)
we obtain
a negatively 
oriented closed arc $A^-_{\varepsilon }$ of
$K''(\delta )$, with starting point $x$, of 
arc length 
$\varepsilon \cdot [ 1 + o(1) ] $,
with $o(\cdot )$ uniform.
Further, analogously to the functions $r_x^+(\cdot )$ and
$s_x^+(\cdot )$ we obtain functions $r_x^-(\cdot )$ and
$s_x^-(\cdot )$. Moreover, there holds the analogue of
(29):
$$
\{ x' \in
K''(\delta ) \mid x' {\text{ follows }} x {\text{ on }}
K''(\delta ) {\text{ in the negative}}
$$
$$
{\text{sense}} \}
\cap B(x, \varepsilon )
=
\big\{ x' \in
K''(\delta ) \mid x' {\text{ follows }} x {\text{ on }}
K''(\delta )
$$
$$
{\text{in the negative sense, and }}
d_{\text{arc}}(x,x') \le \min \{ \delta , s_x^-( 
\varepsilon ) \} \big\} .
\tag 30 
$$

Now
(29)
and
(30)
\thetag{30} 
imply
$$
K''(\delta ) \cap B(x, \varepsilon ) =
\big\{ x' \in
K''(\delta ) \mid x' {\text{ follows }} x {\text{ on }}
K''(\delta ) {\text{ in the positive}}
$$
$$
{\text{sense, and }}
d_{\text{arc}}(x,x') \le \min \{ \delta , s_x^+( 
\varepsilon ) \} \big\} \cup
\big\{ x' \in
K''(\delta )
\mid x' {\text{ follows }} x
$$
$$
{\text{on }}
K''(\delta ) {\text{ in the negative sense, and }}
d_{\text{arc}}(x,x') \le \min \{ \delta , s_x^-( 
\varepsilon ) \} \big\} .
\tag 31 
$$
Applying this to $\varepsilon = \varepsilon (\delta )$
(cf.\
(26))
we get (1) of Lemma 1.3, for ${\text{bd}}\,K$ connected
(cf.\ the beginning of {\bf{1}}) and
$K$ non-compact.

In
(24)
we have seen that
if on $A^+_{\varepsilon } \subset B(x,\varepsilon ) \subset
B(x, 4 \delta )$ a variable
point $x''$ moves from $x$ in the positive orientation, then
$d_X(x,x'')$ is a strictly increasing continuous function of
the position of $x''$. Hence it assumes each value till its
maximum value $\varepsilon $ just once. Analogously we have
that if on $A^-_{\varepsilon }$ a variable
point $x''$ moves from $x$ in the negative orientation, then
$d_X(x,x'')$ is a strictly increasing continuous function of
the position of $x''$. Hence it assumes each value till its
maximum value $\varepsilon $ just once. Applying these for
$\varepsilon := \varepsilon (\delta )$, we get (2) of
Lemma 1.3, for ${\text{bd}}\,K$ connected (cf. the beginning
of {\bf{1}}) and
$K$ non-compact.
As shown in {\bf{5}}, (2) of Lemma 1.3
implies (3) of Lemma 1.3, for ${\text{bd}}\,K$ connected and
$K$ non-compact.
Thus, for ${\text{bd}}\,K$ connected and $K$
non-compact, statements (1), (2), (3) of the lemma are
proved.

\vskip.1cm

{\bf{7.}}
Now we
turn to the case of compact $K$ (then ${\text{bd}}\,K$ is
connected).

Now it makes no sense to say that $x'$ follows $x$ on
$K''(\delta ) = {\text{bd}}\,K$ in the positive, or negative
sense. Therefore we make the following changes in the above
arguments. Choose ${\tilde{x}} \in {\text{bd}}\,K$ such
that it
together with $x$ divides ${\text{bd}}\,K$ to two
arcs of equal lengths. Recall that in
(21)
we have assumed
for the compact case that $\delta \le ({\text{perim}}\,K)/8$.
Then instead of saying
that $x'$ follows $x$ on
$K''(\delta ) = {\text{bd}}\,K$ in the positive, or negative
sense, we will say that $x'$ lies in the positively, or
negatively oriented subarc ${\widehat{x{\tilde{x}}}}$ of
bd\,$K$.
(Then $x,{\tilde{x}}$
lie in both of these arcs, but this causes no problem.)
In the
last paragraph of {\bf{5}} we have shown that the distance
$d(x,x')$ assumes the value $\delta $ for some point of
an arc of bd\,$K$, beginning at $x$ and having arc length
$2\delta $. Observe that by the hypothesis on $\delta $
this arc still lies in one of the subarcs
${\widehat{x{\tilde{x}}}}$ of bd\,$K$. Then we can repeat the
considerations from {\bf{6}}. In
(25)
we will have 
$d_{{\text{arc}}}(x,x') \le
s_x^+(\varepsilon ) = \varepsilon \cdot (1 + o(1)) \le
2\varepsilon \le 2\delta \le ({\text{perim}}\,K)/4$.
We also have the analogous
statement with $s_x^-(\varepsilon )$.

Then with these notational changes all arguments of {\bf{6}}
carry over to the case of compact $K$, and this proves all
statements of this lemma for compact $K$. Then {\bf{6}} and
the previous parts of {\bf{7}} prove (1), (2),
(3) of the
lemma for the case
when ${\text{bd}}\,K$ is connected.

\vskip.1cm

{\bf{8.}}
Second suppose that ${\text{bd}}\,K$ 
has several connected
components. Then ${\text{bd}}\,K$
has at most countably infinitely many
connected components $K'_n$, which are
relatively closed and open subsets of the closed set
${\text{bd}}\,K$. 
Let the component
$K'$ from the statement of this lemma be $K'_1$. 
Then for $n \ge 2$ we have that $\cup _{n \ge 2} K'_n$ is
open and closed in ${\text{bd}}\,K$, hence is closed in $X$.
Now recall that
in a metric space the distance of a compact set, and a closed
set disjoint to it, is positive. This we apply to 
$K''(\delta )$ and $\cup _{n \ge 2} K'_n$.
Then we can choose an 
$\varepsilon _1 > 0$, only depending on $K''(\delta )$ and
$K$, such that the closed $\varepsilon _1$-neighbourhood of 
$K''(\delta )$ is disjoint to $\cup _{n \ge 2} K'_n$.
Thus, in particular,
$$
{\text{for any }} x \in K'' {\text{ we have that }}
B(x, \varepsilon _1) {\text{ is disjoint to }}
\cup _{n \ge 2} K'_n.
\tag 32 
$$

Now let $\varepsilon _2 > 0$ be the value $\varepsilon :=
\varepsilon (\delta )$
which we have
obtained just below
(31),
but with $K$ replaced by the
closed convex set $K^* \supset K$, with boundary $K'$.
Then for
$\varepsilon := \min \{ \varepsilon _1, \varepsilon _2
\} $ statements (1), (2) and (3) of Lemma 1.3 are valid by
{\bf{6}} and {\bf{7}} of this proof and
(32).

{\bf{9.}}
There remains to prove the ``strict monotonicity''
properties of the maps $x \mapsto x_{\varepsilon }^+$ and 
$x \mapsto x_{\varepsilon }^+$, asserted in (4) of the lemma.
We will show this for the map $x \mapsto x_{\varepsilon }^+$
(for the other map the proof is analogous).

By
(24),
for  $x \in K''$ we have $x_{\delta }^+ \in
K''_{2 \delta } \subset K''_{4 \delta } \subset K''(\delta )$.
This implies
$$
(x_{\delta }^+)_{\delta }^+ \in 
K''_{4 \delta } \subset K''(\delta ) .
\tag 33 
$$

By applying
(24)
to $K''(\delta )$, rather than $K''$
(this requires a smaller $\delta $, say, $\delta (0) \in (0,
\delta )$),
on any arc $A$ of 
$K''(\delta )$, of arclength at most $2\delta (0)$, we have the
following. For a proper
subarc $A'$ of $A$, with one endpoint common
with one endpoint of $A$, the corresponding chord length is
strictly smaller than that for $A$. This implies the same
statement for any proper subarc $A'$ of $A$.

Now let
$x^*$ begin to move along $K''(\delta )$, in the
positive sense, from $x \in K''$ till $x_{\delta (0)}^+ \in
K''_{2 \delta }$. Then for
$x^* \ne x$ we cannot have that $(x^*)_{\delta (0)}^+$ belongs 
to the shorter arc ${\widehat{x(x_{\delta (0)}^+)}}$, as follows
from what precedes. Hence $(x^*)_{\delta (0)}^+$, which can
be reached from $x^*$ by moving along $K''(\delta )$
in the positive sense, lies strictly ``beyond''
$x_{\delta (0)}^+$. 

We still have to show that at the beginning of this motion
$(x^*)_{\delta (0)}^+$ cannot move ``too far beyond''
$x_{\delta (0)}^+$.
Let $\beta \in (0, \delta (0)]$ be small. Suppose that,
in the collinear model, $x$ and $x_{\delta (0)}^+$
have images on the negative and positive $\xi $-axis, resp., and
these images are symmetrical w.r.t.\ $0$. Let the length of
the shorter arc ${\widehat{x(x^*)}}$
of $K''(\delta )$ be $\beta $. Then
we assert that $d_X(x_{\delta (0)}^+,(x^*)_{\delta (0)}^+) \le
\beta \left( 1 + o(1) \right) $, for $\beta \to 0$.

Let $\pi $ denote projection on the $\xi $-axis in $X$.

For $X = {\bold{R}}^2, 
H^2$
we have that
$\pi $ is a contraction in $X$. The points $x,x^*,
x_{\delta (0)}^+, $
\newline
$(x^*)_{\delta (0)}^+$ (and
$(x^+_{\delta (0)})^+_{\delta (0)}$, by
(33))
follow each other on
$K''(\delta )$ in the positive order. Therefore,
by the $C^1$-property of $K$, we have that
$x, \pi (x^*), x_{\delta (0)}^+, \pi \left( (x^*)_{\delta (0)}^+
\right) $ follow each other on the $\xi $-axis in the
positive order. 
Then $d_X \left( x, \pi (x^*) \right) \le \beta $ by the
contraction property. This implies $d_X \left( \pi (x^*),
x_{\delta (0)}^+
\right) \ge \delta (0) - \beta $.
Then once more by the
contraction property, we have $d_X \left( \pi (x^*),
\pi \left( (x^*)_{\delta (0)}^+ \right) \right)
\le
d_X \left( x^*, (x^*)_{\delta (0)}^+ \right)
= \delta (0)$.
From these
two inequalities we have 
$d_X \left( 
x_{\delta (0)}^+,
\pi \left( (x^*)_{\delta (0)}^+ \right) \right)
\le \beta
$. Once more by the $C^1$-property of $K$, we have
$d_X \left( x_{\delta (0)}^+,
\left( (x^*)_{\delta (0)}^+ \right) \right)
\le \beta \left( 1 + o(1) \right) $,
for $\beta \to 0$, as asserted.

Next let $X = S^2$.
Since the length of the
shorter arc ${\widehat{x(x^*)}}$ of $K''(\delta )$
is $\beta $, therefore
$d_X \left( x, \pi (x^*) \right) \le d_X(x,x^*) \le
\beta $ (like we had
above).   
Now let $\tilde{x}$ lie on the line $x(x_{\delta (0)}^+)$,
on the
other side of $x_{\delta (0)}^+$, as $x$,
with $d_X(x_{\delta (0)}^+,
\tilde{x}) = \beta $. Then the line passing
through $\tilde{x}$ and orthogonal to the line 
$x(x_{\delta (0)}^+)$ is a great circle $D$, whose distance from
$x$ is
$d_X(x, x_{\delta (0)}^+) + d_X(x_{\delta (0)}^+, \tilde{x}) =
\delta + \beta $.
From above
$d_X(x, x^*) \le \beta $, and then
the distance of $x^*$ from $D$ is at least
$(\delta (0) + \beta ) -  \beta =
\delta (0)$. In particular, by $d_X \left( x^*,
(x^*)_{\delta (0)}^+ \right)
= \delta (0)$, we have that $(x^*)_{\delta (0)}^+ \in B(x^*,
\delta (0))$ cannot lie ``beyond'' $D$. In other words,
$\pi \left( (x^*)_{\delta (0)}^+ \right) $ cannot lie
``beyond'' $D$.
Once more by the $C^1$-property of $K$, we have
$d_X \left( x_{\delta (0)}^+, (x^*)_{\delta (0)}^+
\right) \le \beta \left( 1 + o(1) \right) $,
for $\beta \to 0$, as asserted.
{\bf{QED}}


\vskip.1cm

{\bf{Lemma 1.4.}}
{\it{Assume
{\rm{(1)}}
with $d = 2$.
Let $K$ and $L$ be $C^1$. Let $K'$ and $K''$ be as in Lemma
{\rm{1.3}}.
Let $L'$ and $L''$ be defined
analogously for $L$, as $K'$ and $K''$ 
were defined for $K$ in Lemma {\rm{1.3}}. Then there exists an
$\varepsilon (K,L,K'',L'') > 0$ such that for each
$\varepsilon \in (0, \varepsilon (K,L,K'',L'')]$ 
the following holds.

Let $[x_1,x_2]$ and $[y_1,y_2]$ be chords of $K'$ and $L'$,
resp., of length $\varepsilon $, where $x_2$ follows $x_1$ on
$K'$ in the positive sense, and $y_2$ follows
$y_1$ on $L'$ in the negative sense.
Let at least one of $x_1,x_2$ belong to ${\text{\rm{relint}}}
_{K'} K''$,
and at least one of $y_1,y_2$ belong to ${\text{\rm{relint}}}
_{L'} L''$.
Let us choose $\varphi $ and $ \psi $ 
so, that $\varphi (x_i) = \psi (y_i)$ ($i=1,2$).
Let us consider the shorter arcs ${\widehat{(\varphi x_1)
(\varphi x_2)}}$ and ${\widehat{(\varphi y_2)(\varphi y_1)}}$.
If they are both segments, then this segment is a
subset of $(\varphi K) \cap (\psi L)$. Else
$(\varphi K) \cap (\psi L)$ is the compact convex set with
boundary ${\widehat{(\varphi x_1)(\varphi x_2)}} \cap 
{\widehat{(\varphi y_2)(\varphi y_1)}}$.}}


\vskip.1cm

{\it{Proof.}}
{\bf{1.}}
We begin with the case when ${\text{bd}}\,K$ and
${\text{bd}}\,L$ are connected.

We choose $\varepsilon \in (0, \min \{ \varepsilon (K,K''),
\varepsilon (L,L'') \} ]$
sufficiently small (with $\varepsilon (L,L'')$
defined analogously to $\varepsilon (K,K'')$, cf.\ Lemma 1.3).

In the
second paragraph of {\bf{5}} of the proof of Lemma 1.3 we
have seen the following.
The direction of the oriented segment $S :=
[x,x']$ and the oriented tangent line $l(x')$ of
$K''(\delta )$ at $x'$, both taken in $X$, have in $X$ an
angle, at $x'$, which is $o(1)$, uniformly. We apply this to
$x := x_1,\,\,x' := x_2$, and to $x := x_2,\,\,x' := x_1$,
resp. Analogously, we apply this
to $L''(\delta )$ (defined analogously as $K''(\delta )$,
in {\bf{1}} of the proof of Lemma 1.3),
to $y_1,y_2$ or $y_2,y_1$ in
place of $x,x'$, resp. 
Therefore {\it{the compact convex sets (``half-lens like
domains''), bounded
by $[\varphi (x_1),\varphi (x_2)] =
[\psi (y_1),\psi (y_2)]$,
and by the shorter arcs
${\widehat{(\varphi x_1)(\varphi x_2)}}$
of \,${\text{\rm{bd}}}\,(\varphi K)$ or
${\widehat{(\psi (y_2))(\psi (y_1))}}$ of ${\text{\rm{bd}}}\,
(\psi L)$, resp.,
have at $\varphi (x_1)
= \psi (y_1)$ and $\varphi (x_1) = \psi (y_1)$ angles
which are $o(1)$, uniformly}}. Hence
$$
{\text{the union }} M {\text{ of these two compact convex
sets (``half-lens like domains'')}}
\tag 34 
$$
is also a compact set, which has inner
angles at $\varphi (x_1)
= \psi (y_1)$ and at $\varphi (x_1) = \psi (y_1)$, which are
$o(1)$, uniformly, for $\varepsilon \to 0$.
Hence $M$
is compact convex. (In a suitable orthogonal
coordinate system, in the collinear model, it can be
given as $\{ (\xi , \eta ) \mid \xi \in [\xi _1, \xi _2],\,\,
\eta \in [f(\xi ), g(\xi )]$, with $f$ convex and $g$ concave
on $[\xi_1, \xi _2]$,
and $f(\xi _i) = g(\xi _i) = 0$.)
However, $M$ may have an empty interior.

We are going to prove 
$(\varphi K) \cap (\psi L) = M$, except if both 
${\widehat{(\varphi x_1)(\varphi x_2)}}$ and
${\widehat{(\psi (y_2))(\psi (y_1))}}$ are segments.

\vskip.1cm

{\bf{2.}}
First we prove
$M \subset (\varphi K) \cap
(\psi L)$. It is sufficient to show 
$M \subset \varphi K$. 
$$
{\text{For this it suffices to show that the compact convex
set }} M_1,
$$
$$
{\text{bounded by }} [\varphi x_1, \varphi x_2]
\cup {\widehat{(\psi y_2)(\psi y_1),}} {\text{ is a subset
of }} \varphi K.
\tag 35 
$$

By Lemma 1.3,
${\widehat{(\psi y_2)(\psi y_1)}} \subset B(\psi y_1,
\varepsilon ) = B(\varphi x_1, \varepsilon ) $.
Also by Lemma 1.3,
$\varphi A
:= {\text{bd}}\,(\varphi K) \cap B(\varphi x_1, \varepsilon ) $
is an arc of ${\text{bd}}\,(\varphi K)$,
with endpoints at a
distance $\varepsilon $ from $\varphi x_1$ (one of them being
$\varphi x_2$, the other one is denoted by $\varphi x_0$).
Thus
$$
{\text{no point of }}[{\text{bd}}\,(\varphi K)] \setminus
(\varphi A) = [{\text{bd}}\,(\varphi K)] \setminus
$$
$$
B(\psi x_1, \varepsilon ) {\text{ lies in }}
{\widehat{(\psi y_2) (\psi y_1)}} \subset B(\psi y_1,
\varepsilon ) .
\tag 36 
$$

Also,
$$
{\text{no point of the shorter arc }}
{\widehat{(\varphi x_0)(\varphi x_1)}}
$$
$$
{\text{lies in }} {\widehat{(\psi y_2)(\psi y_1)}},
{\text{ except }} \varphi x_1 = \psi y_1 .
\tag 37 
$$
In fact,
in the collinear model, as a subset of
${\bold{R}}^2$, 
we have
that
${\widehat{(\psi y_2)'(\psi y_1)'}}$ lies close to
$[(\psi y_1)',(\psi y_2)'] =
[(\varphi  x_1)', (\varphi x_2)']$.
(Here $(\cdot )'$ denotes image in the collinear model.)
More exactly, it lies in an angular
domain of vertex $(\varphi x_1)'$, with one leg passing
through $(\varphi x_2)'$, and angular measure $o(1)$.
Similarly, in the collinear model, as a subset of
${\bold{R}}^2$, 
we have that
${\widehat{(\varphi x_0)'(\varphi x_1)'}}$ lies close to
$[(\varphi x_0)', (\varphi x_1)']$. Hence it lies also close to
$[(\varphi x_0^*)', (\varphi x_1)']$, where $\varphi x_0^*$
is the mirror image of $\varphi x_2$ w.r.t.\ $\varphi x_1$.
More exactly, it lies in an angular domain of vertex 
$(\varphi x_1)'$, with one leg passing through
$(\varphi x_0^*)'$, and angular measure
$o(1) + o(1) = o(1)$. Now ${\widehat{(\psi y_2)'(\psi y_1)'}} \cap
{\widehat{(\varphi x_0)'(\varphi x_1)'}}$ lies in
the intersection of these two angular domains, which is,
for $H^2$ and ${\bold{R}}^2$, 
$\{ (\varphi x_1)' \} = \{ (\psi y_1)' \} $. Hence,
for $H^2$ and ${\bold{R}}^2$, 
we have
${\widehat{(\psi y_2)(\psi y_1)}} \cap
{\widehat{(\varphi x_0)(\varphi x_1)}} = \{ \varphi x_1 \}
$, which proves
(37)
in these cases. 
For $S^2$ we still have to take into account that both
${\widehat{(\varphi x_0)(\varphi x_1)}}$ and
${\widehat{(\psi y_2)(\psi y_1)}}$ are subsets of $B(
\varphi x_1, \varepsilon )$. Now the intersection of the
above two angular domains and of $B(
\varphi x_1, \varepsilon )$ is $\{ \varphi x_1 \} $, which
proves
(37)
in this case as well.

If both arcs ${\widehat{(\varphi x_1)(\varphi x_2)}}$
and ${\widehat{(\psi y_2)(\psi y_1)}}$
are equal to the common chord
$[\varphi x_1, \varphi x_2]$
$ = [\psi y_1, \psi y_2]$, then
clearly this common chord is a subset of $(\varphi K) \cap
(\psi L)$, as asserted.
(But this inclusion may be proper.)
This is valid also for ${\text{bd}}\,K$
or ${\text{bd}}\,L$ disconnected.
{\it{Henceforward we exclude the case that both
${\widehat{(\varphi x_1)
(\varphi x_2)}}$ and ${\widehat{(\psi y_2)(\psi y_1)}}$ are
segments (also for the case when the
boundary of $K$ or $L$ is disconnected).}}

If both these
arcs are different from this common chord, then this
common chord strictly
separates these two arcs in $M$,
except for their
endpoints. Now let, e.g., 
${\widehat{(\varphi x_1)(\varphi x_2)}} $
\newline
$ = [\varphi x_1,
\varphi x_2] \ne {\widehat{(\psi y_1)(\psi y_2)}}$. 
Then, in a suitable
rectangular coordinate system in the collinear model,
${\widehat{(\psi y_2)'(\psi y_1)'}} = \{ (\xi, \eta )
\mid \xi \in [\xi _1, \xi _2],\,\,\eta = h(\xi ) \} $. Here
$h$ is concave on $[\xi_1, \xi _2]$
and positive
on $(\xi_1, \xi _2)$, and $h(\xi _i) = 0$. Hence
$\{ (\xi, \eta )
\mid \xi \in [\xi _1, \xi _2],\,\,\eta = h(\xi )/2 \} $
strictly separates
these two arcs in $M$,
except for their endpoints. So, in both of these
cases,
$$
{\text{no point of the shorter arc }}
{\widehat{(\varphi x_1)(\varphi x_2)}}
$$
$$
{\text{lies in }} 
{\widehat{(\psi y_2)(\psi y_1)}}, {\text{ except }}
\varphi x_1 {\text{ and }} \varphi x_2 .
\tag 38 
$$

Then by
(36),
(37)
and
(38),
$$
{\widehat{(\psi y_1)(\psi y_2)}} {\text{ intersects }}
{\text{bd}}\,(\varphi K) {\text{ only at }} \varphi x_1
{\text{ and }} \varphi x_2 .
\tag 39 
$$

\vskip.1cm

{\bf{3.}}
Now we prove $M \subset (\varphi K) \cap (\psi L)$.
By the angular conditions at $\varphi x_1$ and $\varphi x_2$
we see
that in some of their
neighbourhoods, except these points
themselves,
the arc
${\widehat{(\psi y_2)(\psi y_1)}}$ lies in ${\text{int}}\,
(\varphi K)$. Hence by
(39)
the relative interior of this arc lies in
${\text{int}}\,(\varphi K)$. Therefore
the closed
convex curve ${\text{bd}}\,M_1 = 
[(\varphi x_1), (\varphi x_2)] \cup
{\widehat{(\psi y_2)(\psi y_1)}}$ (cf.\
(35)
is a subset of $\varphi K$.
Therefore $M_1 = {\text{conv}}\,({\text{bd}}\,M_1) \subset
\varphi K$.
By
(35)
this proves $M
\subset \varphi K$, thus $M \subset
(\varphi K) \cap (\psi L)$.

\vskip.1cm

{\bf{4.}}
Since in {\bf{2}} we have excluded the case that both
${\widehat{(\varphi x_1)
(\varphi x_2)}}$ and ${\widehat{(\psi y_2)(\psi y_1)}}$ are
segments (also for ${\text{bd}}\,K$ or ${\text{bd}}\,L$
disconnected), therefore $M$ has an
interior point $o$. 
Now we show 
$M \supset (\varphi K) \cap (\psi L)$.
Suppose $z \in 
[(\varphi K) \cap (\psi L)] \setminus M$. Then by {\bf{3}},
$M \subset (\varphi K) \cap (\psi L)$, hence
$o \in {\text{int}}\,M \subset
{\text{int}}\,[(\varphi K) \cap (\psi L)]$.
Then the open
segment $(o,z)$ satisfies
$(o,z) \subset
{\text{int}}\,[(\varphi K) \cap (\psi L)] =
[{\text{int}}\,(\varphi K)] \cap
[{\text{int}}\,(\psi L)]$.
However, it also
intersects ${\text{bd}}\,M$. Suppose, e.g., that
$(o,z)
\cap {\widehat{(\varphi x_1)(\varphi x_2)}} \ne \emptyset $.
Then
$$ 
\emptyset \ne (o,z) \cap
[{\widehat{(\varphi x_1)(\varphi x_2)}}] 
\subset
[{\text{int}}\,(\varphi K)] \cap {\text{bd}}\,(\varphi K) =
\emptyset ,
\tag 40 
$$
a contradiction, showing $M \supset (\varphi K) \cap (\psi L)$.
Hence, by {\bf{3}}, $M = (\varphi K) \cap (\psi L)$.

\vskip.1cm

{\bf{5.}}
Second suppose that ${\text{bd}}\,K$ or ${\text{bd}}\,L$
has several connected
components. This implies that $X$ is ${\bold{R}}^2$ 
or $H^2$.
Recall that in {\bf{2}} we have excluded the case
when both ${\widehat{(\varphi x_1)(\varphi x_2)}}$ and
${\widehat{(\psi y_2)(\psi y_1)}}$ are equal to the common chord
$[\varphi x_1, \varphi x_2] = [\psi y_1, \psi y_2]$, also
for disconnected ${\text{bd}}\,K$ or ${\text{bd}}\,L$.
The sets ${\text{bd}}\,K$
and ${\text{bd}}\,L$ have at most countably infinitely many
connected components $K'_n$ and $L'_m$, resp. Let 
$K' = K'_1$ and 
$L' = L'_1$. 
Then by {\bf{8}} of the proof of Lemma 1.3
$$
{\text{we can choose an }} \varepsilon _1, {\text{ only
depending on }} K''(\delta ),\,\, K,
$$
$$
L''(\delta ) {\text{ and }} L, {\text{ such that the closed }}
\varepsilon _1{\text{-neighbourhood}}
$$
$$
{\text{of }} K''(\delta ) {\text{ is disjoint to }} \cup _{n
\ge 2} K'_n, {\text{ and the closed}}
$$
$$
\varepsilon _1{\text{-neighbourhood of }} L''(\delta )
{\text{ is disjoint to }} \cup _{m \ge 2} L'_m.
\tag 41 
$$

Let $K_n \supset K$ be the closed convex set,
whose boundary is $K'_n$.
Similarly we define $L_m$. Then
$$
(\varphi K) \cap (\psi L) = [\cap _n (\varphi K_n)] \cap
[\cap _m (\psi L_m)] = P \cap Q, {\text{ where}}
$$
$$
P := (\varphi K_1) \cap (\psi L_1) {\text{ and }}
Q := (\cap _{n \ge 2} K_n) \cap (\cap _{m \ge 2} L_m) .
\tag 42 
$$
From the proof of Lemma 1.3 and the previous parts of the
proof of this lemma, for $\varepsilon \in (0, \min \{
\varepsilon (K,K''), \varepsilon (L,L''), \varepsilon _1 \} ]$
sufficiently small, we have
$P \subset B(\varphi x_1, \varepsilon ) \subset Q $.

The {\it{radial function}} of a nonempty
closed convex set $C \subset X$,
w.r.t.\
a point $c \in C$, is defined on
the unit circle of the tangent space of $X$ at $c$, with values
in $[0, \infty ]$, as follows. Its value at some $u$ is the
length of the maximal geodesic segment, starting from $c$, in
direction $u$, contained in $C$. (Recall that now $X$ is
${\bold{R}}^2$ 
or $H^2$!)
The radial
function of $P \cap Q$ w.r.t.\
$\varphi x_1 $ is
the minimum of the radial functions of $P$ and $Q$ w.r.t.\
$\varphi x_1$. The
radial function of $P$ w.r.t.\ $\varphi x_1$
is at most $\varepsilon $,
while that of $Q$ is at least $\varepsilon $.
Hence
the radial function of $P \cap Q$, w.r.t.\
$\varphi x_1$, equals the radial
function of $P$, w.r.t.\
$\varphi x_1$. Hence $(\varphi K) \cap (\psi L) =
P \cap Q = P$.
Thus, using
(42),
the case when ${\text{bd}}\,K$ or
${\text{bd}}\,L$ has
several
connected boundary components,
is reduced to the case of connected
boundaries ($K_1'$ and $L_1'$),
which has been settled in {\bf{4}}.
{\bf{QED}}


\vskip.1cm

{\it{Proof of Theorem {\rm{1}}}}, {\bf{continuation.}}
{\bf{2.}}
Now we continue the proof of the fact, that $(2)
\Longrightarrow (3)$ in Theorem 1.

We will use the notations and hypotheses of Lemma 1.4.
In particular, $[x_1,x_2]$ and $[y_1,y_2]$ have length
$\varepsilon $.

First suppose that not
both shorter arcs ${\widehat{x_1x_2}}$ and
${\widehat{y_1y_2}}$ 
are equal to the corresponding chords. 
Then $(\varphi K) \cap
(\psi L)$ is bounded by the shorter arcs
${\widehat{\varphi (x_1) \varphi (x_2)}}$ and
${\widehat{\psi (y_1) \psi (y_2)}}$, cf.\ Lemma 1.4. These
shorter arcs have lengths at most $\varepsilon \left( 1 + o(1)
\right) $,
cf.\
(17).
Hence
$$
{\text{diam}} \,[ (\varphi K) \cap (\psi L) ] \le
{\text{perim}} \,[ (\varphi K) \cap (\psi L) ] /2 \le
\varepsilon \left( 1 + o(1) \right) .
\tag 43 
$$
Hence $(\varphi K) \cap (\psi L)$
has an arbitrarily small
diameter, for
$\varepsilon \to 0$. Moreover, by
the hypothesis about the arcs, this
intersection has a nonempty interior. Hence, by hypothesis,
it admits a non-trivial congruence. 
Observe that this intersection has just two points of non-smoothness, namely 
$\varphi (x_1) = \psi (y_1)$ and $\varphi (x_2) = \psi (y_2)$. Thus, any
non-trivial congruence admitted by $(\varphi K) \cap
(\psi L)$ is a central symmetry, with centre the midpoint
$o$ of the segment (shorter segment for $S^2$)
joining these two non-smooth points; or it
is an axial symmetry, either
with axis passing
through these two non-smooth
points, or with axis the perpendicular bisector of the
segment with endpoints these two non-smooth points. (For
$S^2$ also $-o$ is a centre of symmetry, but the symmetries
w.r.t.\ the centres $\pm o$ coincide, so we only use $o$.)

Now consider the case that both above arcs
${\widehat{x_1x_2}}$ and ${\widehat{y_1y_2}}$
are equal to the corresponding
chords, which
have length $\varepsilon $. Then $(\varphi K) \cap
(\psi L)$ may strictly contain the common chord $[\varphi x_1,
\varphi x_2] = [\psi y_1,\psi y_2]$,
thus, in particular, its diameter may be not small. 
In this case, therefore, we will
consider, rather than this intersection, this common chord, 
as a degenerate closed convex set (i.e., 
with empty interior). Observe that 
this common chord (in general not equal to $(\varphi K) \cap
(\psi L)$) has an arbitrarily small diameter, and admits all
three non-trivial congruences from the last paragraph.

In both cases, the intersection (in the first case above), 
{\it{or}}\,\,\,the above common chord (in the
second case above), has an arbitrarily small diameter, and 
admits (at least)
one of the three above mentioned non-trivial congruences.

\vskip.1cm


{\bf{Lemma 1.5.}}
{\it{Suppose all hypotheses of Lemma {\rm{1.4}}. Further let $z_1,z_2 \in X$,
such that $d(z_1,z_2) = \varepsilon \in \left( 0, \varepsilon
(K,L,K'',L'') \right) $. We define $\varphi , \psi $
as the unique orientation preserving congruences of $X$ satisfying
$\varphi x_i = \psi y_i = z_i$ for $i = 1,2$.
Then in case {\rm{(A)}}:
${\widehat{(\varphi x_1)(\varphi x_2)}} =
{\widehat{(\psi y_1)(\psi y_2)}} = [z_1,z_2]$, 
we let $M(x_1,x_2,y_1,y_2)$ being this segment.
Else {\rm{(B)}}: let
$M(x_1,x_2,y_1,y_2) := (\varphi K) \cap (\psi L)$. Then the map
$(x_1,x_2,y_1,y_2) \mapsto M(x_1,x_2,y_1,y_2)$ is continuous
from
the set of $(x_1,x_2,y_1,y_2)$'s described in Lemma
{\rm{1.4}} (depending
on $K,K',$
\newline
$K''$ and $L,L',L''$), to the
set of nonempty compact convex sets in $X$, with the
Hausdorff-metric. Equivalently, the map which maps
$(x_1,x_2,y_1,y_2)$ to the image
of $M(x_1,x_2,y_1,y_2)$
in the collinear model,
is continuous in the Hausdorff metric
of
${\bold{R}}^2$ 
(for $S^2$ and ${\bold{R}}^2$), 
and of the
unit circle as a subset of ${\bold{R}}^2$ 
(for $H^2$).}}


\vskip.1cm

{\it{Proof.}}
{\bf{1.}}
$$
{\text{We denote by }} (\cdot )' {\text{ images in the
collinear model, as a subset of }} {\bold{R}}^2 
\tag 44 
$$
(for $S^2$, ${\bold{R}}^2$ 
and $H^2$, resp.).

For $X = H^2$ or  $X = {\bold{R}}^2$  
let $C \subset X$ be compact. For $X = S^2$
let $C$ be a compact set of the open southern hemisphere.
Then by
(6)
and
(7)
we have for distinct
nonempty compact subsets $A,B$ of $C$ the following.
The quotient of the
Hausdorff distances of $A$ and $B$ in $X$, and of 
$A'$ and $B'$ in the collinear model, as a subset of
${\bold{R}}^2$, 
is bounded below and
above. Hence for nonempty compact subsets $A_n$ and $B$ of
$C$,
convergence of $A_n$ to $B$,
in the Hausdorff metric of $X$, is equivalent to the
convergence of $A_n'$ to $B'$, in the Hausdorff metric of the 
collinear model, as a subset of ${\bold{R}}^2$. 
This shows the
equivalence in the last sentence of this lemma.

\vskip.1cm

{\bf{2.}}
We suppose that the midpoint $o$ of $[z_1,z_2]$ is mapped to
$0$ in the collinear model, and that
$z_1'$, or $z_2'$ lies on the negative, or positive
$\xi $-axis, resp. Then the straight line $z'_1z'_2 =
(\varphi x_1)'(\varphi x_2)' = (\psi y_1)'(\psi y_1)'$ is the
{\it{horizontal axis}}, and the perpendicular bisector
straight line of
$[z'_1,z'_2]$ in ${\bold{R}}^2$ 
is the {\it{vertical axis in ${\bold{R}}^2$}}. 
We write $z_i' = (\zeta _i, 0)$.

In both cases (A) and (B) we cut $M(x_1,x_2,y_1,y_2)$ by 
$[z_1,z_2]$ to two nonempty compact convex
parts, namely
$$
M_K, {\text{ bounded by }} {\widehat{(\varphi x_1)
(\varphi x_2)}} {\text{ and }} [z_1,z_2],
$$
$$
{\text{and }} M_L, {\text{ bounded by }} {\widehat{(\psi y_1)
(\psi y_2)}} {\text{ and }} [z_1,z_2] .
\tag 45 
$$
We will
deal with the case of $M_K$ (the case of $M_L$ being
analogous).

By the proof of Lemma 1.3, {\bf{5}}, we have the following.
The angle of
the chord $[z_1,z_2]$, and of the tangents of the arc
${\widehat{(\varphi x_1)(\varphi x_2)}}$, at both endpoints
$z_1$ and $z_2$ are $o(1)$ in $X$. The same statement
holds also
for their respective
images in the collinear
model, as subsets of ${\bold{R}}^2$ 
(by Lemma 1.1).
Hence in the collinear
model, as a subset of ${\bold{R}}^2$, 
we have that
$$
M_K' = \{ (\xi , \eta ) \mid \zeta _1 \le \xi \le \zeta _2,
\,\, f_K(\xi ) \le \eta \le 0 \} . 
\tag 46 
$$
Here $f_K$ is a $C^1$ convex function on $[\zeta _1,
\zeta _2]$,
with
$f_K(\zeta _1) = f_K(\zeta_2) = 0$, and $f_K'(\zeta _1)
= - o(1)$ and $f_K'(\zeta _2) = o(1)$
for $\varepsilon \to 0$ (with
$o(1)$ nonnegative). Then on $[\zeta _1, \zeta _2]$ 
we have that $f$ is non-positive, and $-o(1) \le f' \le o(1)$
for 
$\varepsilon \to 0$. These imply that for
$\xi \in [\zeta _1,\zeta _2]$ we have $f(\xi ) \ge \max \{
-o(1) \cdot (\xi - \zeta _1), o(1) \cdot (\xi - \zeta _2) \}
$. This implies that $M_K'$
lies in the lower semicircle $S'$ of the Thales
circle of $[z_1, z _2]$, meant in ${\bold{R}}^2$ 
(for $\varepsilon $
sufficiently small). Moreover,
$$
M = M_K \cup M_L = {\text{conv}}\,(M_K \cup M_L) ,
{\text{ hence }} M' = {\text{conv}}\,(M_K' \cup M_L') .
\tag 47 
$$
Hence for the {\it{support functions}} $h(M',u)$ etc.\ 
of $M'$, $M_K'$ and $M_L'$
(defined for $u \in S^1$)
we have
$$
h(M',u) = \max \{ h(M_K',u), h(M_L',u) \} .
\tag 48 
$$

We write $C(S^1) =
\{$continuous functions $S^1 \to {\bold{R}} \} $, 
with the
maximum norm, and hence with the topology of uniform
convergence.

Recall that
the Hausdorff distance
of two convex sets in ${\bold{R}}^2$ 
equals
the distance, in the maximum norm, of their support
functions (defined on $S^1$). 
Therefore, in ${\bold{R}}^2$, 
convergence, in the Hausdorff metric,
of nonempty compact convex sets to a limit nonempty compact
convex set,
is equivalent to the uniform
convergence in $C(S^1)$ of their respective support
functions. Now recall that the maximum operation of two
functions in $C(S^1)$
preserves uniform
limits. (I.e., uniform convergence of $f_n$ to $f$ and of 
$g_n$ to $g$ imply uniform convergence of
$\max \{ f_n, g_n \} $ to $\max \{ f, g \} $.) Therefore,
by
(48),
it suffices to prove continuous dependence 
of the support functions
of $M_K'$ and $M_L'$ on $x_1,x_2,y_1,y_2$ in $C(S^1)$.  
Observe that $M_K'$ only depends on $x_1,x_2$, and 
$M_L'$ only depends on $y_1$ and $y_2$. (Since $x_1$ and
$x_2$ are $C^1$ functions of each other, we could even say
dependence only on $x_1$, or only on $x_2$,
whichever lies in ${\text{relint}}_{K'}K''$; 
and similarly for
the $y_i$'s.)

\vskip.1cm

{\bf{3.}}
Therefore we are going to show continuous dependence of
$M_K'$ on $x_1$ and $x_2$. (That of $M_L'$ on $y_1,y_2$ is
shown analogously.)

We distinguish two cases.
\newline
(1): $M_K' \ne [z_1, z_2]$, and
\newline
(2): $M_K' = [z_1, z_2]$.

In case (1) we rewrite convergence of a sequence $A_n 
$ of nonempty compact convex sets to a convex body
$B$ in the Hausdorff metric in 
${\bold{R}}^2$. 
Choose some $b \in {\text{int}}\,B$. Then a
neighbourhood base of $B$ is
$$
\{ B^* \mid \emptyset \ne B^* \subset {\bold{R}}^2 
{\text{ is a
compact convex set, and }} h(B^*,u) \in [h(B,u) - \delta ,
$$
$$
h(B,u) + \delta ] \} =
\{ B^* \mid \emptyset \ne B^* \subset {\bold{R}}^2 
{\text{ is a
compact convex set, and }} h(B^*-b, u)
$$
$$
= h(B^*,u) - \langle b,u \rangle
\in [h(B-b, u) - \delta , h(B-b, u) + \delta ] \} ,
{\text{ for }} \delta \in (0, 1) .
\tag 49 
$$
Since $h(B^* - b,u)$ is bounded below and above, we may
rewrite this as
$$
\{ B^* \mid B^* \subset {\bold{R}}^2 
{\text{ is a nonempty
compact convex set, and }} h(B^*-b, u) \in 
$$
$$
[(1 - \delta )
h(B-b, u), (1 + \delta )h(B-b, u)] \} =
\{ B^* \mid B^* \subset {\bold{R}}^2 
{\text{ is a 
convex}}
$$
$$
{\text{body, and }}
(1 - \delta )(B-b) \subset
(B^*-b) \subset (1 + \delta )(B-b) \} ,
{\text{ for }} \delta \in (0, 1) .
\tag 50 
$$

We denote by $\varrho (B-b,u)$ the radial function of
the
convex body $B-b \subset {\bold{R}}^2$ 
w.r.t.\ $0$, for $0 \in
{\text{int}}\,(B-b)$. (Thus ${\text{bd}}\,(B-b)$
has equation $r = \varrho (B-b,u)$, for $u \in S^1$,
in polar coordinates $u,r$.)
Hence we can rewrite the second half of
(50)
as
$$
\{ B^* \mid B^* \subset {\bold{R}}^2 
{\text{ is a 
convex body, and }} \varrho (B^*-b, u) \in
$$
$$
[(1 - \delta )
\varrho (B-b, u), (1 + \delta )\varrho (B-b, u)] \} ,
{\text{ for }} \delta \in (0, 1) . 
\tag 51 
$$
Since $\varrho (B-b,u)$ is bounded below and above, we can
further rewrite
(51)
as
$$
\{ B^* \mid B^* \subset {\bold{R}}^2 
{\text{ is a 
convex body, and }} \varrho (B^*-b, u) 
$$
$$
\in [\varrho (B-b, u) - \delta ,
\varrho (B-b, u) + \delta ] \} , {\text{ for }}
\delta \in (0, 1)
\tag 52 
$$
or as 
$$
\{ B^* \mid B^* \subset {\bold{R}}^2 
{\text{ is a convex body,
and the distance of }} \varrho (B-b, u) 
$$
$$
{\text{ and }} \varrho (B^*-b, u)
{\text{ in }} C(S^1) {\text{ is at most }} \delta \} ,
{\text{ for }} \delta \in (0, 1) .
\tag 53 
$$

{\bf{4.}}
Let $x_1,x_2$ be as in this lemma, and let $x_1^*,x_2^*$,
with the same properties, be close to $x_1,x_2$.
$$
{\text{Let }} b' \in {\text{int}}\,M_K'.
\tag 54 
$$
Then for $x_i^*$ sufficiently close to $x_i$, we have that
the shorter arc
${\widehat{(\varphi x_1^*)'(\varphi x_2^*)'}}$ of
${\text{bd}}\,(\varphi K)'$
is not equal to
$[(\varphi x_1^*)', (\varphi x_2^*)']$. Moreover,
$b'$ lies on the same
open side of the line $(\varphi x_1^*)'(\varphi x_2^*)'$ as
this shorter arc.

Let $H_-'$
denote the closed
lower halfplane $\eta \le 0$. Let $(H_-')^*$ be
the closed halfplane bounded by the line 
$(\varphi x_1^*)'(\varphi x_2^*)'$, containing the shorter arc
${\widehat{(\varphi x_1^*)'(\varphi x_2^*)'}}$ of
${\text{bd}}\,(\varphi K)'$.
Let $(M_K')^* := (\varphi K)' \cap (H_-')^*$ (a segment of
$(\varphi K)'$). Then, for $x_i^*$ sufficiently close to $x_i$,
we have $b' \in {\text{int}}\,(M_K')^*$. 
Let $(S')^*$ be the
intersection of the Thales circle of 
$[(\varphi x_1^*)', (\varphi x_2^*)']$
(meant in ${\bold{R}}^2$) 
with $(H_-')^*$.
By {\bf{2}}
of this proof we have  $(M_K')^* \subset (S')^* \subset
(H_-')^*$. Hence 
$$
(M_K')^* = (M_K')^* \cap (S')^* = [(\varphi K)' \cap (H_-')^* ]
\cap (S')^*
$$
$$
= (\varphi K)' \cap [(H_-')^* \cap (S')^*] =
(\varphi K)' \cap (S')^* .
\tag 55 
$$
Then for the radial functions we have
$$
\varrho \left( (M_K')^* - b', u \right) = \min \{ 
\varrho \left( (\varphi K)' - b', u \right) ,
\varrho \left( (S')^* - b', u \right) \} . 
\tag 56 
$$

Now observe that, by Lemma 1.3, (4), and elementary geometry,
for $x_i^*$ sufficiently
close to $x_i$, we have 
in the Hausdorff distance of ${\bold{R}}^2$, 
that
$(S')^*$ is sufficiently close to
$S'$. Then, by the equivalence of
(49)
and
(53),
we have that
$\varrho \left( (S')^* - b', u \right) $ is sufficiently close
to $\varrho (S' - b', u)$.
Now recall that the minimum operation of two functions in
$C(S^1)$ preserves uniform limits (in the sense as the maximum
operation of two functions).
Then by
(56),
for $x_i^*$ sufficiently close to $x_i$,
we have that in the $C(S^1)$-norm,
$\varrho \left( (M_K')^* - b', u \right) $
is sufficiently
close to $\varrho (M_K' - b', u)$.
Therefore, once more using the equivalence of
(49)
and
(53),
for $x_i^*$ sufficiently close to $x_i$, we have
in the Hausdorff distance of ${\bold{R}}^2$, 
that $(M_K')^*$ is sufficiently close to $M_K'$.

Then, by {\bf{1}} of this proof,
for $x_i^*$ sufficiently close to $x_i$, we have 
in the Hausdorff distance of $X$ that $M_K^*$
is sufficiently close
to $M_K$ (where $M_K^* \subset X$ is
the inverse image of the subset
$(M_K')^*$ of the model). This proves the statement of the
lemma in case (1) in {\bf{3}}.

{\bf{5.}}
There remains to prove the statement of the
lemma in case (2) in {\bf{3}}. 
Suppose that
$x_i^*$ is sufficiently close to $x_i$ (for given
$\varepsilon = d_X(x_1,x_2) = d_X(x_1^*,x_2^*)$),
and that $x_1^*$
follows $x_1$ on ${\text{bd}}\,(\varphi K)$, e.g., in the
positive sense. Then by Lemma 1.3, (4), the order of our
points on ${\text{bd}}\,(\varphi K)$,
in the positive sense, is $\varphi x_1, \varphi x_1^*,
\varphi x_2, \varphi x_2^*$.
Then
we have to show that the Hausdorff distance, in
${\bold{R}}^2$, 
of $M_K' = [z_1',z_2']$ and $(M_K')^*$ is small.

On the one hand, $M_K$
lies in the $d_X(\varphi x_1, \varphi x_1^*)$-neighbourhood
of $(M_K)^*$, 
where $d_X(\varphi x_1, \varphi x_1^*)$ is small.

So, by {\bf{1}} of this proof,
there remains to show that also $(M_K')^*$ lies in a small
metric
neighbourhood of $M_K'$, meant in ${\bold{R}}^2$, 
containing
the collinear model (${\bold{R}}^2$ 
or the open unit circle).
We are going to prove that
a small metric
neighbourhood of the segment $[(\varphi x_1^*)',
(\varphi x_2)']$, meant in ${\bold{R}}^2$, 
contains $(M_K')^*$,
which implies the above statement. Actually, since the
small metric
neighbourhood of the segment $[(\varphi x_1^*)',
(\varphi x_2)']$, meant in ${\bold{R}}^2$, 
is convex in ${\bold{R}}^2$, 
it suffices to show that it contains
${\text{bd}} \left( (M_K')^* \right) $ (whose convex hull is
$(M_K')^*$).

The boundary ${\text{bd}} \left( (M_K')^* \right) $ consists
of three parts: the segments  $[(\varphi x_1^*)',
(\varphi x_2)']$
and $[(\varphi x_1^*)', (\varphi x_2^*)']$,
and the shorter arc
${\widehat{(\varphi x_2)'(\varphi x_2^*)'}}$. The convex
hull of $[(\varphi x_1^*)',(\varphi x_2)'] \cup 
{\widehat{(\varphi x_2)'(\varphi x_2^*)'}}$ contains 
$[(\varphi x_1^*)', (\varphi x_2^*)']$. Therefore it is sufficient
to show that the small metric
neighbourhood of the segment $[(\varphi x_1^*)',
(\varphi x_2)']$
contains $[(\varphi x_1^*)',(\varphi x_2)']$ and 
\newline
${\widehat{(\varphi x_2)'(\varphi x_2^*)'}}$. The first
containment is obvious, so we only need to show that this 
small metric neighbourhood contains 
${\widehat{(\varphi x_2)'(\varphi x_2^*)'}}$.
In turn, this will be shown if we will show that a small
neighbourhood of $(\varphi x_2)'$ contains 
${\widehat{(\varphi x_2)'(\varphi x_2^*)'}}$.
However, by Lemma 1.3, (4), ${\widehat{(\varphi x_2)'
(\varphi x_2^*)'}}$ has a small Euclidean length.
Since the
chord length is always at most the corresponding arc length, 
therefore the
Euclidean distance of $(\varphi x_2)'$ and any point of
${\widehat{(\varphi x_2)'(\varphi x_2^*)'}}$ is small. 
This proves the statement of the lemma in case (2) in
{\bf{3}}.
{\bf{QED}}


\vskip.1cm

We will call a convex surface in ${\bold{R}}^d$ 
(i.e., the
boundary of a proper closed convex subset of
${\bold{R}}^d$, 
with
non-empty interior)
at some of its points {\it{twice
differentiable}} if the following holds. Locally
it is the graph, in a suitable
rectangular coordinate
system, of a function having a Taylor series expansion of
second degree at this point,
with an error term $o( \| \cdot \| ^2 \| )$. 
By [13],
pp. 31-32 (in both editions), convex surfaces
in ${\bold{R}}^d$ 
are
almost everywhere twice differentiable.
This extends to $S^d$ and $H^d$ by using their collinear
models.


\vskip.1cm

{\bf{Lemma 1.6.}}
{\it{Assume \thetag{1} 
with $d = 2$. Let $K$ and $L$ be $C^1$.
Suppose {\rm{(2)}} of Theorem {\rm{1}}, and all hypotheses of
Lemma {\rm{1.4}}.
Suppose that there exists a sequence $\varepsilon _n \to 0$,
where
each $\varepsilon _n$ is sufficiently small, such that
we have the following. With the notations of Lemmas
{\rm{1.3}} and {\rm{1.4}},
either $K'$, or $L'$ has a chord $[x_1,x_2]$, or $[y_1,y_2]$,
with $x_2$ following $x_1$ in the positive sense,
or $y_2$ following $y_1$ in the negative sense, resp., and
with at least one endpoint in ${\text{\rm{relint}}}_{K'} K''$,
or ${\text{\rm{relint}}}_{L'} L''$, resp.,
such that the following holds.
The chord $[x_1,x_2]$, or $[y_1,y_2]$, is of length
$\varepsilon _n$. Moreover,
the shorter arc determined by this chord, either on $K'$,
or on $L'$,  
is not symmetrical w.r.t.\ the orthogonal bisector
of the chord ({\rm{in particular, the
respective shorter arc is different from
the chord}}). Then this leads to a contradiction.}}


\vskip.1cm

{\it{Proof.}}
{\bf{1.}}
Let $[x_1,x_2]$, or $[y_1,y_2]$ be a chord of $K'$, or of $L'$, with at
least one endpoint in ${\text{relint}}_{K'} K''$, or
${\text{relint}}_{L'} L''$,
and of length
$\varepsilon _n$ (which replaces $\varepsilon $ from
Lemmas 1.4 and 1.5), with $x_2$ following $x_1$ on
${\text{bd}}\,K$ 
in the positive sense,
or $y_2$ following $y_1$ on ${\text{bd}}\,L$
in the negative sense, resp.
{\it{Let $\varphi $ and $\psi $ be chosen so, that}}
$\varphi (x_i)= \psi (y_i) =: z_i$
(for $i=1,2$), where $d_X(z_1,z_2) = \varepsilon _n$.
Then, by Lemma 1.4,
$(\varphi K) \cap (\psi L)$ 
is bounded by the shorter
arcs ${\widehat{\varphi (x_1) \varphi (x_2)}}$ and 
${\widehat{\psi (y_1) \psi (y_2)}}$. 
(Observe that at least one of the arcs ${\widehat{x_1x_2}}$
and 
${\widehat{y_1y_2}}$ is
different from the respective chord. Therefore the case that
$(\varphi K) \cap (\psi L)$
strictly contains this chord, and thus is degenerate, cannot
occur, by Lemma 1.4.)
Then the
intersection $(\varphi K) \cap (\psi L)$ has a nonempty
interior, and by
(43)
it has a diameter at most
$\varepsilon _n \left( 1 + o(1) \right) $, which is
arbitrarily small for $n$ sufficiently large.
Hence it admits some non-trivial
congruence. By the hypothesis of the lemma this
cannot be a symmetry w.r.t.\ the perpendicular bisector
of \,$[z_1,z_2]$, which we call the {\it{vertical axis}}. 
That is, this congruence is a central symmetry w.r.t.\ the
{\it{midpoint $o$ of $[z_1,z_2]$}}, or is an axial symmetry
w.r.t.\ the straight line $z_1z_2$, which we call the
{\it{horizontal axis}} (cf.\ the proof of Theorem 1,
{\bf{2}}).

Observe that both central symmetry w.r.t.\ the
midpoint of $[z_1,z_2]$, and axial symmetry
w.r.t.\ the horizontal axis, cannot occur. Namely, then we
would have also
an axial symmetry w.r.t.\ the vertical axis, which is
excluded by the hypothesis of this lemma.
 
In the case of central symmetry w.r.t.\ $o$,
the two (shorter) arcs ${\widehat{x_1x_2}}$ of
${\text{bd}}\,K$
and ${\widehat{y_1y_2}}$ 
of ${\text{bd}}\,L$, resp., are congruent, with $x_1$
corresponding to $y_2$, and $x_2$ corresponding to $y_1$. 
In case of axial symmetry w.r.t.\ the horizontal
axis, once more the
above arcs are congruent, but now with $x_1$ 
corresponding to $y_1$, and $x_2$ corresponding to $y_2$. 

We will consider the one-sided curvatures, provided
they exist, of $K''$ at $x_i$, in the sense towards $x_{2-i}$
(i.e., of the shorter arc ${\widehat{x_1x_2}}$),
and similarly, of $L''$ at $y_j$, in the sense towards
$y_{2-j}$. Here $x_i \in
{\text{relint}}_{K'} K''$, 
and $y_j \in {\text{relint}}_{L'} L''$.
For both considered symmetries (central, and axial w.r.t.\ the
horizontal axis),
the above considered two
one-sided curvatures exist and are equal
at the corresponding points, or they both
do not exist at the corresponding points. 

Now recall from Lemma 1.4,
that any of $x_1,x_2$, or of $y_1,y_2$ 
could be any point of ${\text{relint}}_{K'} K''$,
or of ${\text{relint}}_{L'} L''$, resp. 

\vskip.1cm

{\bf{2.}}
First suppose the
case that, for all choices of $x_1,x_2,y_1,y_2$, we have
central symmetry w.r.t.\ the midpoint of $[z_1,z_2]$. Then
$\varphi (x_1)$ corresponds by this symmetry to $\psi (y_2)$.
Here
$x_1,y_2$ could be any points of
${\text{relint}}_{K'} K''$ and
${\text{relint}}_{L'} L''$.
Therefore, for
all points of ${\text{relint}}_{K'} K''$ and
${\text{relint}}_{L'} L''$,
the considered one-sided
curvatures exist and are equal, or they do not exist for any points. However,
convex curves -- and surfaces -- are almost everywhere
twice differentiable, in the sense as stated before this
lemma.
This rules out the second case. Now, replacing $x_1,y_2$ by $x_2,y_1$, we
obtain the same for one-sided curvatures, but now in the opposite sense. 
Therefore, at all points of ${\text{relint}}_{K'} K''$
and ${\text{relint}}_{L'} L''$,
the above considered two one-sided curvatures
exist and are equal. Since the curvatures of $K''$ and $L''$
exist almost everywhere, the common values of the two one-sided
curvatures are also equal.

\vskip.1cm

{\bf{3.}}
Second suppose the
case that, for all choices of $x_1,x_2,y_1,y_2$, we have axial symmetry, w.r.t.\ the horizontal axis. Then
$\varphi (x_1)$ corresponds by this symmetry to $\psi (y_1)$. Now,
$x_1,y_1$, and also $x_2,y_2$, could be any points of
${\text{relint}}_{K'} K''$, and
${\text{relint}}_{L'} L''$.
Then, with
this notational change, we repeat the arguments of the preceding paragraph. Thus we 
gain that, at all points of ${\text{relint}}_{K'} K''$
and ${\text{relint}}_{L'} L''$,
the curvatures
exist and are equal.

\vskip.1cm

{\bf{4.}}
As third case, there remains 
the case that, for some choice of $x_1,x_2,y_1,y_2$ we have
central symmetry w.r.t.\ $o$,
and for some other choice of these points
we have axial symmetry w.r.t.\ the horizontal axis. 

$$
{\text{We claim that the configurations of the points }}
x_1,x_2,y_1,y_2 {\text{ in }} K' \times K' 
$$
$$
\times L' \times L' ,{\text{ with }} x_2 {\text{ following }}
x_1 {\text{ in the positive sense, and }} y_2
{\text{ following }} 
$$
$$
y_1 {\text{ in the negative sense, where still we
suppose, that one of }} x_1,x_2 
$$
$$
{\text{belongs to }} {\text{\rm{relint}}}_{K'} K'',
{\text{ and one of }} y_1,y_2 {\text{ belongs to }}
{\text{\rm{relint}}}_{L'} L'', {\text{ and}} 
$$
$$
{\text{that }} d(x_1,x_2)= d(y_1,y_2) = \varepsilon _n,
{\text{ is a connected topological space.}}
\tag 57 
$$
Then
the configuration space of the points 
$x_1,x_2,y_1,y_2$ is the product of the configuration spaces
of the points $x_1,x_2$ and of the points $y_1,y_2$,
and the product of connected spaces is connected. Therefore
$$
{\text{it suffices to show connectedness of the configuration
space of the points }} x_1,x_2
\tag 58 
$$
(then connectedness of the
configuration space of the points $y_1,y_2$ follows
similarly).

Suppose that $x_1$ belongs to ${\text{relint}}_{K'} K''$
(which is by Lemma 1.3
homeomorphic either to an open
segment, or to $S^1$, so it is connected).
Then by Lemma 1.3, (4), $x_2$ depends
continuously on $x_1$, hence also the ordered pair $(x_1,x_2)$
depends continuously on $x_1$. Since a continuous image of a
connected space is connected, therefore
these ordered pairs form a connected space.
We have the analogous statement if
$x_2 \in {\text{relint}}_{K'} K''$.
Moreover, these two
connected spaces intersect, provided some $x_1,x_2 \in
{\text{relint}}_{K'} K''$
have a distance
$\varepsilon _n$. This happens if the arclength of 
$K''$ in $X$ is greater than $2 \varepsilon _n$
(cf.\
(24)),
which can be supposed.
Now it suffices to recall that the union of two intersecting
connected spaces is itself connected. This ends the proof
of
(58),
hence of
(57).

Further, we claim that
$$
{\text{the set of configurations of the points }} x_1,x_2,y_1,
y_2, {\text{ for which one of the}}
$$
$$
{\text{considered symmetry properties holds, is a closed
subset of }} K' \times K' \times L' \times L' .
\tag 59 
$$
In fact, by Lemma 1.5, the map $(x_1,x_2,y_1,y_2) \mapsto
M(x_1,x_2,y_1,y_2)'$ (this is the image,
in the collinear model,
of $M(x_1,x_2,y_1,y_2)$) is continuous in the Hausdorff
metric of ${\bold{R}}^2$ 
(for $X = S^2, {\bold{R}}^2$),  
or of
the unit circle, as a subset of ${\bold{R}}^2$ 
(for $X = H^2$).
We suppose that $z'_1$ and $z'_2$ lie in the negative and
positive $\xi $-axis, resp., and
$(z'_1 + z'_2)/2 = 0$. Then $o' = 0$, and the above horizontal
and vertical axes (cf.\ {\bf{1}} of this proof)
are mapped in the collinear model into the
$\xi $- and $\eta $-axes in ${\bold{R}}^2$, 
resp. 

Then central symmetry of $M(x_1,x_2,y_1,y_2)$ w.r.t.\ $o$,
i.e., central symmetry of
$M(x_1,x_2,$
\newline
$y_1,y_2)'$ w.r.t.\ $0$, can
be expressed via the support function of $M(x_1,x_2,y_1,$
$y_2)'$
as
$h \left( M(x_1, \right.$
\newline
$\left. x_2,y_1,y_2)', u \right) =
h \left( M(x_1,x_2,y_1,y_2)', -u \right) $.
Analogously, symmetry of $M(x_1,x_2,y_1,y_2)$ w.r.t.\ 
the horizontal axis, i.e., symmetry of
$M(x_1,x_2,y_1,y_2)'$
w.r.t.\ the $\xi $-axis, can be expressed as
$h \left( M(x_1,x_2,y_1,y_2)', (u_1, u_2) \right) = 
h \left( M(x_1,x_2, 
y_1,y_2)', (u_1, -u_2) \right) $. Clearly
both of these properties are preserved by (uniform) convergence
of the functions $(x_1,x_2,y_1,y_2) \mapsto
h \left( M(x_1,x_2,y_1,y_2)', u \right) $
to a limit function. This proves
(59).

Further, the union of these two nonempty closed subsets is the
entire space of all above configurations of the points $x_1,x_2,y_1,y_2$.
By
(57),
these two closed subsets must intersect.
That is, we must
have a configuration, that simultaneously has both the central symmetry w.r.t.\ $o$,
and the axial symmetry w.r.t.\ the horizontal axis. This, however, 
contradicts the second paragraph of {\bf{1}} of
the proof of this lemma.

\vskip.1cm

{\bf{5.}}
So the third case (investigated in {\bf{4}})
cannot occur. Therefore we must have
either the first, or the second
case (investigated in {\bf{2}} and {\bf{3}}, resp.).
Both had the conclusion that, 
at all points of ${\text{relint}}_{K'} K''$ and
${\text{relint}}_{L'} L''$,
the curvatures
exist and are equal.
In other words, both $K''$ and
$L''$ have equal constant curvatures, i.e., both are arcs of congruent cycles
(including entire compact cycles, i.e., circles),
or are segments.

Recall that $K''$, or $L''$ were arbitrary compact subarcs of $K'$, or $L'$, if
$K'$, or $L'$ were homeomorphic to ${\bold{R}}$, 
and they were
equal to
$K'={\text{bd}}\,K$, or $L'={\text{bd}}\,L$, if $K'$, or $L'$ 
was homeomorphic to 
$S^1$. Thus in both cases, $K'$ and $L'$ are congruent cycles, or
are straight lines. However, this contradicts the assumptions
about ``not axial symmetry of some shorter arc w.r.t.\ the
orthogonal bisector line of the corresponding chord''
of this lemma. 
{\bf{QED}}


\vskip.1cm

{\bf{Lemma 1.7}}
{\it{Assume \thetag{1} 
with $d = 2$.
Let $K$ be $C^1$. Let $K'$ and $K''$ be as
in Lemma {\rm{1.3}}. 
Suppose that for each sufficiently small $\varepsilon > 0$
we have the following. 
For any
chords $[x_1,x_2]$ of $K'$,
with $x_2$ following $x_1$ on $K'$ in the positive sense,
and
with at least one endpoint in
${\text{\rm{relint}}}_{K'} K''$,
and with length of these chords being
$\varepsilon $, the following holds.
The shorter arcs determined by these chords, on $K'$,
are symmetrical w.r.t.\ the perpendicular
halving straight line of the chord ({\rm{in particular, the
respective shorter arcs may coincide with
these chords}}). Then $K'$ is a cycle, or a straight line.
Moreover, if $K'$ is a circle or paracycle, then $K$ is
a circle (disk) or a paracircle, resp.
For $L'$ and $L$ there holds the analogous
statements, with $y_2$ following $y_1$ on $L'$ in the negative
sense.}}


\vskip.1cm

{\it{Proof.}}
Observe that the perpendicular
halving straight line of such a chord also 
halves the shorter respective arc, and is perpendicular to it
at its midpoint (it cannot touch this arc).

By the statement before Lemma 1.6 convex curves are almost
everywhere twice differentiable, in the sense given there,
hence have curvatures almost everywhere. 
Now let $x' \in {\text{relint}}_{K'} K''$,
such that $K'$ is twice differentiable at
$x'$. Further, let $x'' \in {\text{relint}}_{K'} K''$
be arbitrary.
Then there exist $x'=x_1, \dots , x_n=x''
\in {\text{relint}}_{K'} K''$,
following each other in the same
sense, and
such, that the distance of $x_i$ and 
$x_{i+1}$ (in $X$)
is less than $\varepsilon $, for $i=1, \dots , n-1$. 
Then $x_i$ and 
$x_{i+1}$ are symmetrical w.r.t.\
the perpendicular bisector of the chord $[x_i,x_{i+1}]$.
Then $x_i$ and
$x_{i+1}$ are symmetrical also w.r.t.\ the
perpendicular bisector of some other 
chord. For this the corresponding shorter arc $I'$
contains the closed shorter arc $I={\widehat{x_ix_{i+1}}}$ 
in its relative interior w.r.t.\ $K'$, with
$I'$ being only slightly longer
than $I$. Moreover, these
two arcs have the same midpoint, and the chord corresponding
to $I'$ is still shorter than $\varepsilon $.
Further, also $I'$ is symmdtrical w.r.t.\ the perpendicular
bisector of $[x_i, x_{i + 1}]$.

Then by induction
one sees that $K''$ is twice differentiable at each $x_i$,
and the curvatures of $K''$ at each $x_i$ are equal to the
same number $\kappa  \ge 0$. In particular,  
the curvature of $K'$ at any $x'' \in {\text{relint}}
_{K'} K''$
equals the constant $\kappa $.
Recall from Lemma 1.3
that $K''$ was equal to $K'={\text{bd}}\,K$, if $K$ was
compact and then $K'$ was
homeomorphic to $S^1$. Moreover, $K''$ was
any compact arc of $K'$ if $K$ was
noncompact and then $K'$ was
homeomorphic to ${\bold{R}}$. 
Hence the curvature of $K'$ at any of its points
equals the constant $\kappa $.

For $K$ noncompact (thus 
$X = {\bold{R}}^2, 
H^2$) we have that
$K'$ joins two infinite points 
(possibly coinciding). Both for $K$ compact and noncompact,
and for any $X$, we have that
$K'$ is a maximal
$C^2$ curve of constant curvature, i.e., a cycle,
or a straight line, as asserted. If $K'$ is a circle or a
paracycle, then ${\text{conv}}\,K' \subset K$. If here we
have an equality, then $K$ is a circle or a paracircle.
Otherwise,
there is a $k \in K \setminus ({\text{conv}}\,K')$. Then for
$k^* \in {\text{int}}\,({\text{conv}}\,K')$ we have that the
open segment $(k^*,k)$, which lies in ${\text{int}}\,K$,
intersects $K'$, which lies in ${\text{bd}}\,K$. This is a
contradiction. (Also cf.\ the
last paragraph of the proof of Lemma 1.9
in [8].)
This ends the proof of the lemma.
{\bf{QED}}


\vskip.1cm

{\it{Proof of Theorem {\rm{1}}}}, {\bf{continuation.}}
{\bf{3.}}
Observe that Lemmas 1.6 and 1.7, with their complementary
hypotheses, together prove 
$(2) \Longrightarrow (3)$ in Theorem 1.

The particular case of $(2) \Longrightarrow (3)$, for
central symmetry, now follows easily. If the connected
components $K'$ and $L'$ are not congruent, one of them
is strictly convex. Now let
$\varphi K$ and $\psi L$ touch each other, so that 
$K'$ and $L'$ touch each other. Pushing them slightly
towards each other, the new intersection has a small diameter,
hence is centrally symmetric. It has two points of
non-smooothness, the common endpoints of the boundary arcs on 
$K'$ and $L'$. The central symmetry should exchange
these two common endpoints, hence also these two boundary 
arcs. However, these two boundary 
arcs have different curvatures, so this is impossible.
This contradiction implies the  
equality of the curvatures of $K'$ and $L'$.


\vskip.1cm

Next we turn to the investigation of the implication
$(3) \Longrightarrow (1)$ in Theorem 1, under the respective
hypotheses.


\vskip.1cm


{\bf{Lemma 1.8.}}
{\it{Assume
{\rm{(1)}}
with $d = 2$, and {\rm{(3)}} of Theorem
{\rm{1}}. For $X = H^2$,
let both for $K$ and $L$, the infimum of the
positive curvatures of its connected
boundary components be positive,
and let it have at most one boundary component with $0$
curvature. Then for $X = S^2,\,\, {\bold{R}}^2, 
\,\, H^2$,
{\rm{(1)}} of Theorem {\rm{1}} holds.}}


\vskip.1cm

{\it{Proof.}}
{\bf{1.}}
For $X=S^2$, $(3)$ of Theorem 1 clearly implies $(1)$
of Theorem 1. (Observe that for both $K$ and $L$ halfspheres,
$(1)$ of Theorem 1 holds vacuously.)

\vskip.1cm

{\bf{2.}}
For $X = {\bold{R}}^2$, 
supposing (3) of Theorem 1, we have 
the following. The closed convex sets in
${\bold{R}}^2$, 
whose boundaries are disconnected, are just the
parallel strips.
Furthermore, the closed 
convex sets in ${\bold{R}}^2$, 
with connected boundaries, 
whose boundaries are cycles or straight lines, are just
circles or
half-planes, resp. Thus, any of $K$ and $L$ can be a
circle, a
parallel strip, or a half-plane.
If one of $K$ and $L$ is a circle, 
then $(\varphi K) \cap (\psi L)$ is axially symmetric, hence
(1) of Theorem 1 holds.
If each of $K$ and $L$ is either a parallel strip or a
half-plane, then there does not exist  
$(\varphi K) \cap (\psi L)$, with nonempty interior, and of
arbitrarily small diameter. Hence now (1) of Theorem 1 holds
vacuously.

\vskip.1cm

{\bf{3.}}
Let $X=H^2$, and suppose (3) of Theorem 1. Then
each of $K$ and $L$ is either a circle, or a paracircle, or has boundary components
which are hypercycles or straight lines. 
The infimum of the positive curvatures of its boundary
components is of course the same infimum, taken only for the
hypercycle
components.
By the hypothesis of the lemma, 
the distances, for which these hypercycles are distance lines, have an
infimum $c>0$, say. Moreover, there may be still at most one straight 
line component, both for $K$ and for $L$. 

Let, e.g., $K_1$ and $K_2$ be two boundary components of
$K$, and let $x_1 \in K_1$, and $x_2 \in K_2$. Let
$K' \subset K$ and $L' \subset L$ be defined,
as the nonempty closed convex sets (possibly with empty interiors), 
bounded by all the straight lines
for which the boundary components are distance lines, and by the at most one
straight line component. In particular, $K_1$ and
$K_2$ are distance lines for $K_1'$ and $K_2'$, with a
non-negative signed distance.
Then the segment $[x_1,x_2]$ intersects both $K_1'$ and $K_2'$,
at
points $x_1' \in K'_1$ and $
x_2' \in K'_2$, and for the distances we have $d(x_1,x_2)
\ge d(x_1,x_1') + d(x_2',x_2)
\ge c$. This means that the distances of the different
boundary components both of $K$,
and of $L$, are bounded from below by $c$. The same holds vacuously for circles
and paracircles. Hence, if diam\,$[ (\varphi K) \cap
(\psi L) ] < c$, then 
$(\varphi K) \cap (\psi L)$ is compact, and is 
bounded by portions of only one boundary
component of $\varphi K$, and of $\psi L$. 

Thus $(\varphi K) \cap (\psi L)$ is the intersection of
two sets, both being a circle, a paracircle, or a convex domain bounded by a
hypercycle, including a half-plane. 
Recall that a circle, and a paracircle are axially symmetric
w.r.t.\ any straight line passing through their centres. Thus, 
if both above sets are a circle or a 
paracircle, then their intersection is axially symmetric. There
remain the cases when one set is a convex set bounded by a hypercycle, and the
other one is a circle, a paracircle, or a convex set bounded by a hypercycle.  
In the first case an axis of symmetry of the
intersection is the straight line passing through the centre of the circle, and
orthogonal to the base line of the hypercycle. In the second case, by
compactness of the intersection, the centre of the paracircle cannot lie at an
endpoint of the base line. Therefore an axis of symmetry of the
intersection is the straight line passing through the centre of the paracircle, 
and orthogonal to the base line of the hypercycle. In the third case, again by
compactness of the intersection, having interior points,
the base lines of the
hypercycles are 
ultraparallel. Moreover, the hypercycles lie on those closed sides of their
base lines, as the other base line. Therefore, the unique
straight line orthogonal to both base lines
is an axis of symmetry of the intersection.
{\bf{QED}}


\vskip.1cm

Last we turn to the investigation of $(3) \not\Longrightarrow
(2)$ in Theorem 1, under the respective hypotheses. 


\vskip.1cm

{\bf{Lemma 1.9.}}
{\it{Assume
(1)
with $d = 2$, and {\rm{(3)}} of Theorem
{\rm{1}} for
$K$. Let $X = H^2$, and let $K$ be not a circle or a
paracircle.
Let us prescribe in any way the curvatures of the 
hypercycle and straight line connected boundary components
of $K$ (with
multiplicity),
so that
the infimum of the positive
curvatures is $0$, or there are two $0$ curvatures. Then there
exists a $K$, with these prescribed curvatures of its
hypercycle and straight line
boundary components (with
multiplicity), and an $L$, such that {\rm{(3)}} of Theorem
{\rm{1}}
holds for $K$ and $L$, but {\rm{(2)}} of Theorem {\rm{1}}
does not
hold for them.}}


\vskip.1cm

{\it{Proof.}}
{\bf{1.}}
We begin with an example $K = K_0$ and $L = L_0$,
where each of bd\,$K_0$ and bd\,$L_0$
consists of two straight lines.
We consider the collinear model. Let 
$k_1,k_2 \subset H^2$ be distinct parallel straight 
lines, with axis of symmetry $k$. Let, for $i = 1,2$, the
points
$x_i,y_i \in k_i$ be symmetric w.r.t.\ $k$, with all
six pairwise
distances at most $\varepsilon $. Moreover, let 
the point $x_i$ separate $y_i$ and
the common infinite point of $k_1$ and $k_2$.
Then $x_1x_2y_2y_1$ is a
symmetrical quadrangle of diameter at most $\varepsilon $
(since $X = H^2$). Moreover, it 
is the intersection of the closed convex sets $K_0$,
bounded by $k_1$ and $k_2$, and $L_0$, bounded by
the straight lines $x_1x_2$ and $y_1y_2$.

Let us consider a small generic
perturbation $x_1'x_2'y_2'y_1'$ 
of the quadrangle
$x_1x_2y_2y_1$, where $x_1'$ is the small
perturbation of $x_1$ etc., satisfying
$x_i',y_i' \in k_i$. Then by this perturbation $K_0$ goes
over to
the closed convex set $K_0'$ bounded by the parallel but
distinct lines $x_1'y_1' =
k_1$ and $x_2'y_2' = k_2$, i.e., $K_0' = K_0$.
Moreover, $L_0$
goes over to
the closed convex set $L_0'$
bounded by the ultraparallel straight
lines $x_1'x_2'$ and $y_1'y_2'$.
Then any non-trivial congruence admitted by the perturbed
quadrangle preserves both pairs of opposite sides
(separately the parallel and the ultraparallel ones),
and preserves
the above separation properties. Therefore $x_1'$ (and $y_1'$)
has as image
either itself, or $x_2'$ (and $y_2'$).
Then the congruence is an identity, or it
exchanges $x_1',x_2'$ as well as $y_1',y_2'$. The second case
is
only possible if $d(x_1',y_1') = d(x_2',y_2')$.
Generically this equality
does not hold, so generically 
a non-trivial congruence admitted by the perturbed
quadrangle does not exist. Let us fix some such generic
quadrangle $x_1'x_2'y_2'y_1'$, which admits
no non-trivial congruence.
We may suppose that
diam\,$(x_1'x_2'y_2'y_1') \le 2 \varepsilon $.

\vskip.1cm

{\bf{2.}}
Suppose that the set (with multiplicity) of the positive 
curvatures of the connected hypercycle boundary
components $K_i$ of $K$ is prescribed, and has infimum $c=0$,
or there are at
least two $0$ curvatures. Then we make the
following generalization of the above example. We begin with
constructing $K$. {\it{These
hypercycles $K_i$ are distance
lines, with base lines $K_i^*$, and for prescribed
distances $c_i$}} (cf.\ \S 3). We consider a closed convex set
$K'$, 
bounded in the
collinear model by the prescribed many, at least two,
but at most countably infinitely many chords $K_i^*$
of the collinear model
circle, one for each $i$.
Let, with at most one exception, these chords occur in
disjoint pairs having exactly
one common endpoint. Hence ${\text{int}}\,K' \ne \emptyset $.
Then we replace these chords $K_i^*$ by the corresponding
distance lines $K_i$, outwards from $K'$. If
\newline
(1)
there are two $0$ curvatures, then
the corresponding chords $K_i^*$ should 
occur in an above
pair, and if
\newline
(2) $c=0$, then there should be above pairs of 
$(K_i^*)$'s, for which both distances $c_i$ 
are arbitrarily small.
\newline
We define $L$ as $L_0'$, and we define $\psi $ as identity.

In case (1) there are two boundary components $K_i$, with
$c_i = 0$, hence satisfying $K_i = K_i^*$, with a common
infinite point.
Recall that any three distinct points of the boundary circle
of the model (collinear or conformal) can be taken over to any
other three distinct points of the boundary circle, of the
same orientation, by (the extension of)
some orientation-preserving congruence. Therefore we may
choose $\varphi $ so that it takes these two $K_i$'s
to the above $k_1,k_2$. Also we may
suppose that, in the collinear model, the images of
$k_1$ and $k_2$
enclose a small angle. Moreover, the image of
$\psi L = \psi L_0' = L_0'$ lies in a small neighbourhood
(meant in ${\bold{R}}^2$, 
which contains the collinear
model circle)
of the common infinite point of $k_1$ and $k_2$.
Then the images of all other $K_j$'s in the collinear model
lie far from the common infinite point of $k_1$ and $k_2$,
therefore $(\varphi K) \cap (\psi L)$ remains unchanged if we
delete all these $K_j$'s from ${\text{bd}}\,K$.
Then $(\varphi K) \cap (\psi L)$
equals the above fixed (generic)
quadrangle $x_1'x_2'y_2'y_1'$, which admits no
non-trivial congruence. This proves the lemma for case (1).

In case (2) there are two boundary components $K_i$, with
both respective $c_i$'s arbitrarily small, with a common
infinite point. Consider the respective base lines $K_i^*$,
and let us choose $\varphi $ so that it takes these two
$(K_i^*)$'s to the above $k_1,k_2$. Like in case (1), we
may suppose
that, in the collinear model, $k_1$ and $k_2$ enclose a small
angle. Moreover,
$\psi L_0 = L_0'$ lies in a small neighbourhood
(in ${\bold{R}}^2$, 
containing the collinear model)
of the common infinite point of $k_1$ and $k_2$.
Like in case (1), deletion of all other $K_j$'s from
${\text{bd}}\,K$ lets $(\varphi K) \cap (\psi L)$ unchanged.
Then $(\varphi K) \cap
(\psi L)$ is an arc-quadrangle. It is bounded by two
hypercycle arcs lying on the $K_i$'s,
each very close to $k_1$ and to $k_2$, resp., and
by two segments lying on the lines $x_1'x_2'$ and $y_1'y_2'$.
Then $(\varphi K) \cap (\psi L)$
is a very small perturbation of the quadrangle 
$x_1'x_2'y_2'y_1'$.
Hence we may suppose that diam\,$[ (\varphi K)
\cap (\psi L) ] \le 3 \varepsilon $. Moreover, the quadrangle 
$x_1'x_2'y_2'y_1'$ admits no non-trivial congruences.

Suppose that for both respective
$c_i$'s anyhow small $(\varphi K) \cap (\psi L)$
admits a non-trivial congruence. Then in limit
(of some subsequence) we
would obtain a non-trivial congruence admitted by the
quadrangle
$x_1'x_2'y_2'y_1'$, contrary to the choice of this quadrangle.
Hence, for both $c_i$'s sufficiently small,
itself $(\varphi K) \cap (\psi L)$ cannot admit any
non-trivial congruences. This proves the lemma for case (2),
and the proof of the lemma is finished.
{\bf{QED}}


\vskip.1cm

{\it{Proof of Theorem}} 1, {\bf{continuation.}}
{\bf{4.}}
Now the previous parts of the
proof of Theorem 1, and Lemmas 1.2 and 1.6-1.9
prove all statements of Theorem 1.
{\bf{QED}}


\vskip.1cm

{\it{Proof of Theorem {\rm{2}}.}}
For the first three statements 
we have the evident implications $(3) \Longrightarrow (1)
\Longrightarrow (2)$. Now we show the
remaining implication $(2) \Longrightarrow (3)$.
Let (2) hold.
Then by Theorem 1, $(2) \Longrightarrow (3)$,
and compactness of $K$ and $L$ we have that
$K$ and $L$ are circles. If both $K$ and $L$ are halfspheres,
then (3) holds. Else we may apply (2), yielding that the
circles $K$ and $L$ are congruent.
That is, (3) holds.

For the last five statements
we have the evident implications $(8) \Longrightarrow (4)$
and $(4) \Longrightarrow (5) 
\Longrightarrow (7)$, and $(4) \Longrightarrow (6) 
\Longrightarrow (7)$. The
remaining implication $(7) \Longrightarrow (8)$ follows from
Theorem 1, $(2) \Longrightarrow (3)$.
{\bf{QED}}


\vskip.1cm

{\it{Proof of Theorem {\rm{3}}.}}
The implication $(2) \Longrightarrow (1)$ is evident. 

For the implication $(1) \Longrightarrow (2)$ we
observe that evidently (1) of Theorem 3 implies (2) of
Theorem 1, and then we can apply
Theorem 1, $(2) \Longrightarrow (3)$. In {\bf{2}} 
of the proof of Lemma 1.8, for $X = {\bold{R}}^2$, 
we have seen that (3) of Theorem 1 implies that
any of $K$ and $L$ can be a circle, a parallel strip, 
or a half-plane. Moreover, if one of $K$ and $L$ is a circle,
then $(\varphi K) \cap (\psi L)$ is axially symmetric.
For the remaining cases observe that the intersection of two
parallel strips is always centrally symmetric, and the
intersection of two half-planes is always axially symmetric.
However, the intersection of a parallel strip and a
half-plane, with nonempty interior,
admits in general no non-trivial congruence.
Thus $(1) \Longrightarrow (2)$ holds.

The two particular cases, with central, or axial
symmetries in $(1)$, follow by easy discussions.
{\bf{QED}}


\vskip.1cm

{\it{Proof of Theorem {\rm{4}}.}}
{\bf{1.}}
The implication $(2) \Longrightarrow
(1)$ is evident, so we turn to the proof of $(1)
\Longrightarrow (2)$.

\vskip.1cm

{\bf{2.}}
Observe that (1) of Theorem 4 implies (2) of Theorem 1, and
(2) of Theorem 1 implies, by Theorem 1, (3) of Theorem 1.
By Theorem 1, for the case of central symmetries in
(2) of Theorem 1, the connected components of the boundaries
of $K$ and $L$ are congruent.
%
That is, 
{\it{$K$ and $L$ are either two congruent circles, or
two paracircles, 
or all their boundary components are either
congruent hypercycles,
or straight lines}}. However, in the case of straight lines,
their total number
is finite, by the hypothesis of the theorem.

The case 
that $K$ and $L$ are paracircles is clearly impossible.
Namely, we may choose
$\varphi $ and $\psi $ so, that $\varphi K = \psi L$, and
then their 
intersection is
a paracircle. However, this has exactly one point at infinity,
hence is not
centrally symmetric.

\vskip.1cm

In the next two lemmas we are going to show that also the
case
of (finitely many)
straight lines, and the case of congruent
hypercycles is impossible.


\vskip.1cm

{\bf{Lemma 4.1}}
{\it{Assume
{\rm{(1)}}
with $d = 2$ and let $X = H^2$. 
Then the case, when all connected components of the boundaries
both of $K$ and $L$ are straight lines,
when, by hypothesis, their total number is finite,
is impossible.}}


\vskip.1cm

{\it{Proof.}}
Now it will be convenient to use the collinear model for $H^2$. Then, in this
model, both $K$ and $L$ are bounded by finitely many
non-intersecting open
chords of the boundary circle of the model. Possibly we have chords with
common end-points. Let $K_1$, or $L_1$ be some connected component of 
${\text{bd}}\,K$, or ${\text{bd}}\,L$, resp. We may choose $\varphi $
and $\psi $ so, that $\varphi K_1 = \psi L_1 = (\varphi K) \cap (\psi L)$, 
and this line
contains the centre of the model. Thus
$\varphi K$ and $\psi L$ lie on the
opposite sides of this straight
line. Let us change
$\varphi $ and $\psi $ a bit, so that in the model
$\varphi K$ and $\psi L$ rotate a 
little bit about the centre of the model.
(Suppose in Remark 1 that $\varphi K_1$ is the vertical axis,
and $\varphi K$ or $\psi L$ lies on the right, or left hand
side of $\varphi K_1$, resp. Then in 
case (1) from
Remark 1 we rotate $\varphi K$ in the negative sense and
$\psi L$ in the positive sense, while in case (2) from
Remark 1 conversely.)
We will not use new notations for the new orientation preserving congruences,
but will retain the old ones $\varphi $ and $\psi $.

Let the intersection $C$ of the closed
half-circles of the collinear model circle,
bounded by $\varphi K_1$, or $\psi L_1$,
and containing $\varphi K$, or $\psi L$, in their new positions, resp,
satisfy the following. It does not contain any end-point of any chord, which
in the model represents some boundary component of $\varphi K$ or $\psi L$, 
except of
course one end-point of $\varphi K_1$, and one end-point of
$\psi L_1$.
By the finiteness hypothesis,
this can be
attained, and implies the following. The set $C$
does not intersect the closure in ${\bold{R}}^2$ 
of
any
other boundary  components of $\varphi K$, or of
$\psi L$ (i.e.,
different from $\varphi K_1$, or $\psi L_1$, resp.) 
than those, which satisfy the following properties (1) and
(2). 
\newline
(1) They are in
the collinear model chords of the model circle with one common end-point
with the chords
$\varphi K_1$, or $\psi L_1$, resp. Moreover, this/these common
end-point/s lie in $C$ (i.e., is/are endpoint/s of the circular arc
corresponding to $C$). 
\newline
(2) From this/these connected 
component/s of the boundaries
only a/ half-line/s is/are in $C$.

Then $(\varphi K) \cap (\psi L)$ is, in
the collinear
model, either
\newline
(a) a sector of the model circle, or
\newline
(b) a triangle, with two sides parallel,
and having two finite vertices, or
\newline
(c) a quadrangle, with
opposite sides parallel.  
\newline
Case (a) gives a set having
exactly one non-smooth boundary point. If it were centrally symmetric, this
boundary point would be the centre of symmetry, which is a contradiction. In
case (b) 
we have a set having exactly one point at infinity, hence it is not centrally
symmetric, which is a contradiction.
In case (c), if there were a centre of symmetry, 
that would be an inner point of our set. Then
one side and its centrally symmetric image side would span
ultraparallel straight 
lines. However, the lines spanned by
any two sides of this quadrangle are
either intersecting, 
or parallel. So we have a contradiction in each of the three cases.

This ends the proof of the lemma.
{\bf{QED}}


\vskip.1cm

{\bf{Lemma 4.2}}
{\it{Assume
(1)
with $d = 2$ and let $X = H^2$.
Then the
case, when all connected components of the boundaries   
both of $K$ and $L$ are congruent hypercycles (degeneration
to straight lines not admitted), is impossible.}}

\vskip.1cm


\vskip.1cm

{\it{Proof.}}
Denote by
$l > 0$ the common value of the distance, for which these
boundary component hypercycles are
distance lines for their base lines.

Again, it will be convenient to consider the collinear model.
Both for $K$ and $L$, remove from it the union of the
convex hulls of its boundary components, thus
obtaining the closed
convex sets $K_0 \,\,(\subset K)$ and $L_0\,\,(\subset L)$.
If there are two boundary
components of $K$ or $L$ with both infinite points common,
then $K_0$ or $L_0$ is a straight line, which we consider as
doubly counted. Else $K_0$ or $L_0$ is a closed convex set
with interior points, and its boundary components are the
base lines of the boundary components of $K$ or $L$, resp.

The parallel domain of $K_0$, or $L_0$, with
distance $l$,
contains $K$, or $L$, resp. However, also
these parallel domains are 
contained in $K$, or $L$, 
resp. Namely, if, e.g., $z \in K_0$, then also $z \in K$.
If however the
distance of a point $z \not\in K_0$ 
from some $x \in K_0$ is at most $l$, then the segment $[z,x]$
intersects some boundary component ${\tilde{K}}_{0,1}$
of $K_0$, say, at a point
$x'$. Then $d(z, x') \le d(z,x) \le l$,
hence the distance
of $z$ from its own projection to ${\tilde{K}}_{0,1}$ is also
at most
$l$. Therefore $z$ lies (not strictly)
between ${\tilde{K}}_{0,1}$ and the
respective
boundary component ${\tilde{K}}_1$ of $K$, hence
$z \in {\text{cl\,conv}}\,{\tilde{K}}_1 \subset K$. (If there
are two
such ${\tilde{K}}_1$'s, then the above statement holds for one
of them.)
That is, we have (in both cases)
$$
{\text{the parallel domain of }} K_0 , {\text{ or }} L_0 ,
{\text{ with distance }}
l, {\text{ equals }} K, {\text{ or }}
L, {\text{ resp.}}
\tag 60 
$$

Let $\varphi K_1$, or $\psi L_1$ denote a boundary component
of $\varphi K$, or $\psi L$, whose base line is denoted by
$\varphi K_{0,1}$, or $\psi L_{0,1}$, resp. Let
$\varphi K_1^*$ or $\psi L_1^*$ denote the closed convex set
bounded by $\varphi K_1$ or $\psi L_1$, resp.
Let us suppose that $\varphi K$ and $\psi L$ 
are in such a position, that $\varphi K_1^*$ and $\psi L_1^*$
have exactly one infinite point in common (which must be the
unique common infinite point of $\varphi K_1$ and
$\psi L_1$). This can be attained by applying some
orientation preserving congruences $\varphi , \psi $.

$$
{\text{Let }} M := (\varphi K_1^*) \cap (\psi L_1^*).
\tag 61 
$$
By the conformal model,
{\it{the set $M$ is bounded by some arcs of $\varphi K_1$ and
$\psi L_1$,
having one common infinite endpoint, which is the only
infinite point of $M$}}.
Evidently $\varphi K \subset \varphi K_1^*$, and $\psi L
\subset \psi L_1^*$, hence
$$
(\varphi K) \cap (\psi L) \subset M .
\tag 62 
$$

We are
going to show that also $(\varphi K) \cap (\psi L)
\supset M$.
It will suffice to show $M \subset \varphi K$
(the other inclusion is proved analogously).

The straight
line $\varphi K_{0,1}$ cuts $H^2$ into two closed half-planes
$\varphi H_{K,0,1}^{\pm }$, with $\varphi K_1 \subset
\varphi H_{K,0,1}^+$. 

For $z \in M \cap (\varphi H_{K,0,1}^+)$ we have $z \in
(\varphi K_1^*) \cap (\varphi H_{K,0,1}^+) =
{\text{cl\,conv}}\,
(\varphi K_1) \subset \varphi K$, as stated.

For $z \in M \cap (\varphi H_{K,0,1}^-)$ we have $z \in
(\psi L_1^*) \cap (\varphi H_{K,0,1}^-)$. Now the straight
line $\psi L_{0,1}$ cuts $H^2$ into two closed half-planes
$\psi H_{L,0,1}^{\pm }$, with $\psi L_1 \subset
\varphi H_{L,0,1}^+$.
Then $\varphi H_{K,0,1}^- \subset 
\psi H_{L,0,1}^+$, hence $z \in (\psi L_1^*) \cap 
(\psi H_{L,0,1}^+) = {\text{cl\,conv}}\,(\psi L_1)$.
Therefore for the projection $\psi q$ of $z$ to
$\psi L_{0,1}$ we have $d(z,\psi q) \le l$.
By $z \in \varphi H_{K,0,1}^-$
we have that $z$ and $\psi q$ are (not strictly)
separated by $\varphi
K_{0,1}$, hence the segment $[z, \psi q]$ intersects $\varphi
K_{0,1}$, at a point $\varphi p$, say. 
Then
$$
{\text{dist}}\,(z,\varphi K_0) \le
{\text{dist}}\,(z, \varphi K_{0,1}) \le d(z, \varphi p) \le
d(z, \psi q) \le l ,
\tag 63 
$$
hence, by
(60),
$z \in \varphi K$, as stated.

Therefore $z \in M$ implies $z \in \varphi K$, i.e., $M
\subset \varphi K$. Similarly $M \subset \psi L$, 
thus $M \subset (\varphi K) \cap (\psi L)$. Then,
by
(62),
$$
(\varphi K) \cap (\psi L) = M .
\tag 64 
$$
As written above (just below
(61)),
the set $M$ has just one point at infinity, which implies that
it cannot be centrally symmetric.
This ends the proof of the lemma.
{\bf{QED}}


\vskip.1cm

{\it{Proof of Theorem}} 4, {\bf{continuation.}}
{\bf{3.}}
Now Theorem 4 follows from the previous parts of the
proof of Theorem 4, and from
Lemmas 4.1-4.2.
{\bf{QED}}


\vskip.1cm

{\bf{Acknowledgements.}} 
The authors express their gratitude to I. B\'ar\'any, for
carrying the problem, and bringing the two authors
together; to V. Soltan, for having sent to the second named
author the manuscript of [15],
prior to its publication; and to K. B\"or\"oczky
(Sr.), for calling the attention of the second named author 
to the fact that, in Theorem
1, our original hypothesis $C^3_+$ was unnecessary. Following
K. B\"or\"oczky's
arguments, the
authors finally succeeded to eliminate all smoothness
hypotheses. 



\vskip.6cm

{\centerline{REFERENCES}}

\vskip.6cm







[1] D. V. Alekseevskij, E. B. Vinberg, A. S. Solodovnikov,
{\it{Geometry of spaces of constant curvature}},
Geometry II (Ed. E. B. Vinberg), Encyclopaedia
Math. Sci. {\bf{29}}, 1-138, Springer, 
Berlin
1993,
MR {\bf{95b:}}{\rm{53042b}}

[2] R. Baldus,
{\it{Nichteuklidische Geometrie, Hyperbolische Geometrie
der Ebene,}} {\rm{4-te Aufl., Bearbeitet und erg\"anzt von
F. L\"obell
(Non-Euclidean geometry, hyperbolic geometry of the
plane, 4th ed., revised and expanded by F. L\"obell,
in German),}}
Sammlung G\"oschen (G\"oschen Collection)
970/970a, de Gruyter,
Berlin, 
1964,
MR {\bf{29\#}}{\rm{3936}}

[3] T. Bonnesen, W. Fenchel,
{\it{Theorie der konvexen K\"orper,}} {\rm{Berichtigter
Reprint (Theory of convex bodies, corrected reprint, in
German),}}
Springer,
Berlin-New York,
1974,
MR {\bf{49\#}}{\rm{9736}}

[4] R. Bonola,
{\it{Non-Euclidean geometry, a critical and
historical study of its developments,}} {\rm{Translation with
additional appendices by H. S. Carslaw.
With a Supplement containing the
G. B. Halstead translations of}} {\it{``The science of absolute space''}}
{\rm{by 
J. Bolyai and}} {\it{``The theory of parallels''}} {\rm{by N. Lobachevski}}
Dover Publs. Inc.,
New York, N.Y.,
1955,
MR {\bf{16-}}{\rm{1145}}

[5] H. S. M. Coxeter
{\it{Non-Euclidean Geometry}}, 6th ed., 
Spectrum Series, The Math. Ass. of America,
Washington, DC,
1998,
MR {\bf{99c:}}{\rm{51002}}

[6] E. Heil, H. Martini,
{\it{Special convex bodies,}} 
In: Handbook of Convex
Geometry (eds. P. M. Gruber, J. M. Wills), North-Holland,
Amsterdam etc., 1993, Ch. 1.11,
347-385, 
MR {\bf{94h:}}{\rm{52001}}
 
[7] R. High,
{\it{Characterization of a disc, Solution to problem 1360 (posed by
P. R. Scott),}}
Math. Magazine
{\bf{64}} 
(1991),
353-354 

[8] J. Jer\'onimo-Castro, E. Makai, Jr.,
{\it{Ball characterizations in spaces of constant curvature,}}
Studia Sci. Math. Hungar.
{\bf{55}}
(2018),
421-478, 
MR {\bf{3895882}}{\rm{}} 

[9] J. Jer\'onimo-Castro, E. Makai, Jr.,
{\it{Ball characterizations in planes and spaces of constant
curvature, II}},
manuscript in preparation

[10] J. Jer\'onimo-Castro, E. Makai, Jr.,
{\it{Ball characterizations in Euclidean spaces}},
manuscript in preparation 


[11] H. Liebmann
{\it{Nichteuklidische Geometrie}}, {\rm{3-te Auflage
(Non-Euclidean geometry, 3rd ed., in German)}}
de Gruyter,
Berlin,
1923,
Jahr\-buch Fortschr. Math. {\bf{49}}, 390

[12] O. Perron,
{\it{Nichteuklidische Elementargeometrie der Ebene}}
{\rm{(Non-Euclidean elementary geometry of the plane, in
German)}},
{\rm{Math. Leitf\"aden}}
Teubner,
Stuttgart,
1962
MR {\bf{25\#}}{\rm{2489}}

[13] R. Schneider,
{\it{Convex bodies: the Brunn-Minkowski theory; 
Convex bodies: the Brunn-Minkowski theory, Second expanded
edition}}, {\rm{Encyclopedia of
Math. and its Appls., Vol.}} {\bf{44}}; {\bf{151}},
Cambridge Univ. Press,
Cambridge,
1993; 2014,
MR {\bf{94d:}}{\rm{52007}}; {\bf{3155183}}

[14] V. Soltan,
{\it{Pairs of convex bodies with
centrally symmetric intersections of translates}},
Discrete Comput. Geom. 
{\bf{33}}
(2005),
605-616, 
MR {\bf{2005k:}}{\rm{52012}} 

[15] V. Soltan,
{\it{Line-free convex bodies
with centrally symmetric intersections of translates}},
Revue Roumaine Math. Pures Appl. 
{\bf{51}}
(2006)
111-123, 
MR {\bf{2007k:}}
{\rm{52010}}.
{\rm{Also in: Papers on Convexity and Discrete geometry, Ded. to T. 
Zamfirescu on the occasion of his
60th birthday, Editura Academiei Rom\^ane, Bucure\c sti, 2006, 411-423}}

[16] J. J. Stoker,
{\it{Differential Geometry}},
New York Univ., Inst. Math. Sci.,
New York,
1956

[17] I. Vermes,
{\it{\"Uber die synthetische Behandlung der
Kr\"ummung und des Schmieg\-zy\-kels der ebenen Kurven in der
Bolyai-Lobatschefskyschen Geometrie}} {\rm{(On the synthetic
treatment of the curvature and the osculating cycle in the
Bolyai-Lobachevski\u\i an geometry, in German)}}
Studia Sci. Math. Hungar.
{\bf{28}}
(1993),
289-297, 
MR {\bf{95e:}}{\rm{51030}} 

[18]
{\it{Beltrami-Klein model}},
Wikipedia

[19]
{\it{Hyperbolic triangle}},
Wikipedia

[20]
{\it{Poincar\'e disk model}},
Wikipedia




\end